\documentclass[aos,nameyear,seceqn,dvips]{arximspdf}
\usepackage{graphicx}

\doi{10.1214/10-AOS835}
\volume{38}
\issue{6}
\pubyear{2010}
\firstpage{3811}
\lastpage{3836}

\makeatletter

\newtheorem{theorem}{Theorem}
\newproclaim{remark}{Remark}
\newproclaim{example}{Example}[section]

\newcommand{\T}{{\mathrm{T}}}

\newcommand{\hbbeta}{\hat{\bbeta}}
\newcommand{\ve}{\varepsilon}
\newcommand{\bA}{\mathbf A}
\newcommand{\bX}{\mathbf X}
\newcommand{\bfX}{\mathbf X}
\newcommand{\bfZ}{\mathbf Z}
\newcommand{\bv}{\mathbf v}
\newcommand{\bR}{\mathbf{R}}
\newcommand{\bzeta}{{\bolds{\zeta}}}
\newcommand{\balpha}{{\bolds{\alpha}}}
\newcommand{\bbeta}{{\bolds{\beta}}}
\newcommand{\bdelta}{{\bolds{\delta}}}
\newcommand{\Sba}{\Sigma_{\mathop{\balpha}_{1}}}
\newcommand{\Sbb}{\Sigma_{\mathop{\bbeta}_{1}}}
\newcommand{\wideti}{\widetilde}
\newcommand{\wideh}{\hat}

\newcommand{\bfD}{\mathbf D}

\newcommand{\GXX}{\Gamma_{\wideti X_1\wideti X_1}}
\newcommand{\GXZ}{\Gamma_{\wideti X_1\wideti Z_1}}
\newcommand{\GZZ}{\Gamma_{\wideti Z_1\wideti Z_1}}

\newcommand{\BIC}{\operatorname{BIC}}
\newcommand{\df}{\operatorname{DF}}
\newcommand{\SE}{\operatorname{SE}}
\newcommand{\mse}{\operatorname{MSE}}
\newcommand{\rss}{\operatorname{RSS}}
\newcommand{\sgn}{\operatorname{sgn}}
\newcommand{\diag}{\operatorname{diag}}

\makeatother

\begin{document}
\begin{frontmatter}

\title{Estimation and testing for partially linear single-index models}
\runtitle{Partially linear single-index models}

\begin{aug}
\author[A]{\fnms{Hua} \snm{Liang}\thanksref{t1}\corref{}\ead[label=e1]{hliang@bst.rochester.edu}},
\author[A]{\fnms{Xiang} \snm{Liu}\thanksref{t2}\ead[label=e2]{xliu@bst.rochester.edu}},
\author[B]{\fnms{Runze} \snm{Li}\thanksref{t3}\ead[label=e3]{rli@stat.psu.edu}}
\and
\author[C]{\fnms{Chih-Ling} \snm{Tsai}\ead[label=e4]{cltsai@ucdavis.edu}}
\thankstext{t1}{Supported in part by NIH/NIAID Grant AI59773 and NSF
Grant DMS-08-06097.}
\thankstext{t2}{Supported in part by the Merck Quantitative Sciences
Fellowship Program.}
\thankstext{t3}{Supported by National Science Foundation Grants DMS-03-48869
and
NIDA, NIH Grants R21 DA024260 and P50 DA10075. The content is solely
the responsibility of
the authors and does not necessarily represent the official views
of the NIDA or the NIH.}
\runauthor{Liang, Liu, Li and Tsai}
\affiliation{University of Rochester, University of Rochester, Pennsylvania State University
and~University~of~California, Davis}

\address[A]{H. Liang\\
X. Liu\\
Department of Biostatistics\\
\quad and Computational Biology\\
University of Rochester\\
Rochester, New York 14642\\
USA\\
\printead{e1}\\
\phantom{E-mail: }\printead*{e2}} 
\address[B]{R. Li\\
Department of Statistics\\
Pennsylvania State University\\
University Park, Pennsylvania 16802\\
USA\\
\printead{e3}}
\address[C]{C.-L. Tsai\\
Graduate School of Management\\
University of California, Davis\\
Davis, California 95616\\
USA\\
\printead{e4}}
\end{aug}

\received{\smonth{10} \syear{2009}}
\revised{\smonth{5} \syear{2010}}

%
\begin{abstract}
In partially linear single-index models, we obtain the
semiparametrically efficient profile least-squares estimators of
regression coefficients. We also employ the smoothly clipped absolute
deviation penalty (SCAD) approach to simultaneously select variables
and estimate regression coefficients. We show that the resulting SCAD
estimators are consistent and possess the oracle property.
Subsequently, we demonstrate that a proposed tuning parameter selector,
BIC, identifies the true model consistently. Finally, we develop a
linear hypothesis test for the parametric coefficients and a
goodness-of-fit test for the nonparametric component, respectively.
Monte Carlo studies are also presented.
\end{abstract}

%


\end{frontmatter}

\section{Introduction}\label{sec:intr}

Regression analysis is commonly used to explore the relationship
between a response variable $Y$ and its covariates $Z$. For the sake of
convenience, one often employs a linear regression model $E(Y|Z) =Z^{\mathrm{T}}\balpha$
to assess the impact of the covariates on the response,
where $\balpha$ is an unknown vector and $E$ stands for
``expectation.'' In practice, however, the linear assumption may not be
valid. Hence, it is natural to consider a single-index model $E(Y|Z)
=\eta(Z^{\T}\balpha)$, in which the link function $\eta$ is unknown.
Accordingly, various parameter estimators of single-index models have
been proposed [e.g., Powell, Stock and Stocker (\citeyear{powstosto1989}); Duan and Li
(\citeyear{duali1991}); H{\"a}rdle, Hall and Ichimura (\citeyear{harhalich1993}); Ichimura (\citeyear{ich1993}); Horowitz and
H{\"a}rdle (\citeyear{horhar1996}); Liang and Wang (\citeyear{liawan2005})]. Detailed discussion and
illustration of the usefulness of this model can be found in Horowitz
(\citeyear{hor1998}). Although the single-index model plays an important role in data
analysis, it may not be sufficient to explain the variation of
responses via covariates $Z$. Therefore, Carroll et al.
(\citeyear{carfangijwan1997}) augmented the single-index model in a linear form with
additional covariates $X$, which yields a partially linear single-index
model (PLSIM), $E(Y|Z,X)=\eta(Z^{\T}\balpha)+X^{\T}\bbeta$. When
$Z$ is
scalar and $\balpha=1$, the PLSIM reduces to the partially linear
model, $E(Y|Z,X)=\eta(Z)+X^{\T}\bbeta$ [see Speckman (\citeyear{spe1988})]. A
comprehensive review of partially linear models can be found in
H\"ardle, Liang and Gao (\citeyear{harliagao2000}).

To estimate parameters in partially linear single-index models, Carroll
et al.~(\citeyear{carfangijwan1997}) proposed the backfitting algorithm. However, the
resulting estimators may be unstable [see Yu and Ruppert (\citeyear{yurup2002})] and
undersmoothing the nonparametric function is necessary to reduce the
bias of the parametric estimators. Accordingly, Yu and Ruppert (\citeyear{yurup2002})
proposed the penalized spline estimation procedure, while Xia and
H{\"a}rdle (\citeyear{xiahar2006}) applied the minimum average variance estimation
(MAVE) method, which was originally introduced by Xia et al.
(\citeyear{xiatonlizhu2002}) for dimension reduction. Although Yu and Ruppert's procedure is
useful, it may not yield efficient estimators; that is, the asymptotic
covariance of their estimators does not reach the semiparametric
efficiency bound [Carroll et al. (\citeyear{carfangijwan1997})]. In addition, these estimators
need to be solved via an iterative procedure; that is, iteratively
estimating the nonparametric component and the parametric component. We
therefore propose the profile least-squares approach, which obtains
efficient estimators and provides the efficient bound. Moreover, the
resulting estimators can be found without using the iterative procedure
mentioned above and hence may reduce the computational burden.

In data analysis, the true model is often unknown; this allows the
possibility of selecting an underfitted (or overfitted) model, leading
to biased (or inefficient) estimators and predictions. To address this
problem, Tibshirani (\citeyear{tib1996}) introduced the least absolute shrinkage and
selection operator (LASSO) to shrink the estimated coefficients of
superfluous variables to zero in linear regression models.
Subsequently, Fan and Li (\citeyear{fanli2001}) proposed the smoothly clipped absolute
deviation (SCAD) approach that not only selects important variables
consistently, but also produces parameter estimators as efficiently as
if the true model were known, a property not possessed by the LASSO.
Because it is not a simple matter to formulate the penalized function
via the MAVE's procedure, we employ the profile least-squares approach
to obtain the SCAD estimators for both parameter vectors $\balpha$ and
$\bbeta$. Furthermore, we establish the asymptotic results of the SCAD
estimators, which include the consistency and oracle properties [see
Fan and Li (\citeyear{fanli2001})]. Simulation results are consistent with theoretical
findings.

After estimating unknown parameters, it is natural to construct
hypothesis tests to assess the appropriateness of the linear
constraint hypothesis as well as the linearity of the nonparametric
function. We demonstrate that the resulting test statistics are
asymptotically chi-square distributed under the null hypothesis. In
addition, simulation studies indicate that the test statistics perform
well. The rest of this paper is organized as follows. Section~\ref{sec:prof}
introduces the profile least-squares estimators and the penalized SCAD
estimators. The asymptotic properties of these estimators are obtained.
Section~\ref{sec:test} presents hypothesis tests and their large-sample properties.
Monte Carlo studies are presented in Section~\ref{sec:nume}. Section~\ref{sec:disc} concludes the
article with a brief discussion. All detailed proofs are deferred to
the \hyperref[sec:appe]{Appendix}.

\section{Profile least-squares procedure}\label{sec:prof}

Suppose that $\{(Y_i, Z_i, X_i), i=1,\ldots, n\}$ is a random sample
generated from the PLSIM
%
\begin{equation}\label{eq:2.1}
Y=\eta(Z^{\T}\balpha)+X^{\T}\bbeta+\ve,
\end{equation}
where $Z$ and $X$ are $p$-dimensional and
$q$-dimensional covariate vectors, respectively,
$\balpha=(\alpha_1,\ldots,\alpha_p)^{\T}$,
$\bbeta=(\beta_1,\ldots,\beta_q)^{\T}$, $\eta(\cdot)$ is an unknown
differentiable function, $\ve$ is the random error with mean zero and
finite variance $\sigma^2$, and $(Z^{\T},X^{\T})^{\T}$ and $\ve$ are
independent. Furthermore, we assume that $\|\balpha\|=1$ and that the
first element of $\balpha$ is positive, to ensure identifiability. We
then employ the profile least-squares procedure to obtain efficient
estimators and SCAD estimators in the following two subsections,
respectively.

\subsection{Profile least-squares estimator}\label{sec:plse}

In semiparametric models, Severini and Wong (\citeyear{sevwon1992}) applied the profile
likelihood approach to estimate the parametric component. We adapt this
approach to estimate unknown parameters of partially linear
single-index models. To begin, we re-express model (\ref{eq:2.1}) as
\[
Y_i^*=\eta(\Lambda_i) + \ve_i,
\]
where $Y^*_i=Y_i-X_i^\T\bbeta$ and $\Lambda_i=Z_i^\T\balpha$.
Then, for
a given $\bzeta=(\balpha^\T, \bbeta^\T)^\T$, we employ the local linear
regression technique [Fan and Gijbels (\citeyear{fangij1996})] to estimate $\eta$, that
is, to minimize
%
\begin{equation}\label{eq:2.2}
\sum_{i=1}^n\{a+b(\Lambda_i-u)+X_i^\T\bbeta-Y_i\}^2K_h(\Lambda_i-u)
\end{equation}
with respect to $a$, $b$, where
$K_h(\cdot)=1/hK(\cdot/h)$, $K(\cdot)$ is a kernel function and $h$ is
a bandwidth.

Let $(\hat{a},\hat{b})$ be the minimizer of (\ref{eq:2.2}). Then,
%
\begin{equation}\label{eq:2.3}
\hat{\eta}(u,\bzeta) =\hat{a}= \frac{K_{20}(u,\bzeta
)K_{01}(u,\bzeta) -
K_{10}(u,\bzeta)K_{11}(u,\bzeta)}
{K_{00}(u,\bzeta)K_{20}(u,\bzeta)-K_{10}^2(u,\bzeta)},
\end{equation}
where $K_{jl}(u,\bzeta) = \sum_{i=1}^n K_h(Z_i^{\T}\balpha-u)
(Z_i^{\T}\balpha-u)^j(X_i^{\T}\bbeta-Y_i)^l$ for $j=0,1,2$ and $l=0,1$.
Subsequently, following Jennrich's (\citeyear{jen1969}) assumption (a), there exists
a profile least-squares estimator $\hat{\bzeta}=(\hat{\balpha}{}^{\T},
\hat{\bbeta}{}^{\T})^{\T}$ obtained by minimizing the following profile
least-squares function:
%
\begin{equation}\label{eq:2.4}
Q(\bzeta)=\sum_{i=1}^n\{Y_i-\hat{\eta}(Z_i^{\T}\balpha,
\bzeta)-X_i^{\T}\bbeta\}^2.
\end{equation}
It is noteworthy that
we apply the Newton--Raphson iterative method to find the estimator
$\hat{\bzeta}$; this technique is distinct from the more commonly used
iterative procedure in partially single-index models which iteratively
updates estimates of the nonparametric component and the parametric
component obtained from their corresponding objective functions. In
addition, our proposed profile least-squares approach allows us to
directly introduce the penalized function, as given in the next
subsection.

To study the large-sample properties of parameter estimators, we
consider the true model with an unknown parameter vector
$\bzeta_0=(\balpha_0^\T,\bbeta_0^\T)^\T$. In addition, we assume that
$\balpha_0^\T$ and $\bbeta_0^\T$ have the same dimensions as their
corresponding parameter vectors $\balpha^\T$ and $\bbeta^\T$ in the
candidate model of this subsection. Moreover, we introduce the
following notation: $A^{\otimes2}=AA^\T$ for a matrix $A$,
$\Lambda=Z^{\T}\balpha$, $\wideh{\Lambda}=Z^{\T}\wideh{\balpha}$,
$\Lambda_0=Z^\T\balpha_0$, $\wideti{\xi}=\xi-E(\xi|\Lambda)$,
$\wideti{\xi}_0=\xi-E(\xi|\Lambda_0)$ and $\wideh{\xi}=\xi-\wideh {E}(\xi
|\Lambda)$ for
any random variable (or vector) $\xi$, where $\wideh E(\xi|\Lambda)$ is
the local linear\vspace*{1pt} estimator of $E(\xi|\Lambda)$. For example, $\wideti Z=Z-
E(Z|\Lambda)$ and $\wideh X=X-\wideh E(X|\Lambda)$. We next present six
conditions and then obtain the weak consistency and asymptotic
normality of the profile least-squares estimators.
\begin{longlist}
\item The
function $\eta(\cdot)$ is differentiable and not constant on the
support ${\mathcal U}$ of $Z^\T\balpha$.
\item The function
$\eta(z^\T\balpha)$ and the density function of $Z^{\T}\balpha$,
$f_{\balpha}(z)$, are both three times continuously differentiable with
respect to $z$. The third derivatives are uniformly Lipschitz
continuous over ${\mathcal A}\subset{\mathcal R}^{p}$ for all
$u\in\{u=z^\T\balpha\dvtx \balpha\in{\mathcal A}, z\in{\mathcal Z}\subset{\mathcal R}^{p}\}$.
\item $E(|Y|^{m_1})<\infty$ for some $m_1\ge3$. The conditional
variance of $Y$ given $(X,Z)$ is bounded and bounded away from $0$.
\item The kernel function $K(\cdot)$ is twice continuously
differentiable with the support $(-1, 1)$. In addition, its second
derivative is Lipschitz continuous. Moreover, $\int u^jK(u)\,du=1$ if
$j=0$ and $=0$ if $j=1$.
\item The bandwidth $h$ satisfies
$\{nh^{3+3/{(m_1-1)}}\}\log^{-1}n\to\infty$ and $nh^8\to0$ as
$n\to\infty$.
\item $E(\wideti X^{\otimes2})$ and $E\{\wideti
Z\eta'(\Lambda)\}^{\otimes2}$ are positive-definite, where $\eta
'$ is
the first derivative of $\eta$.
\end{longlist}

\begin{theorem}\label{th:th1}
Under the regularity conditions \textup{(i)--(vi)}, with probability tending to
one, $\hat{\bzeta}$ is a consistent estimator of $\bzeta$. Furthermore,
$\sqrt{n}(\hat{\bzeta}-\bzeta_{0}) \to N(0, \sigma^2 \mathbf{D}^{-1})$ in
distribution, where $\mathbf{D}=E[\{\eta'(\Lambda_0)\wideti Z_0^\T, \wideti
X_0^{\T}\}^{\T}]^{\otimes2}$. Moreover, $\hat{\bzeta}$ is a
semiparametrically efficient estimator.
\end{theorem}

Because $(Z^{\T},X^{\T})^{\T}$ is independent of $\ve$, we are able to
demonstrate that the asymptotic variance of Xia and H{\"a}rdle's
(\citeyear{xiahar2006}) minimum average variance estimator is $\sigma^2\mathbf{D}^{-1}$,
which indicates that MAVE is also an efficient estimator.

After having estimated $\balpha$ and $\bbeta$, we obtain the following
estimator of~$\eta(u)$:
\[
\hat{\eta}(u)=\hat{\eta}(u,\hat{\bzeta}) =
\frac{K_{20}(u,\hat{\bzeta})K_{01}(u,\hat{\bzeta}) -
K_{10}(u,\hat{\bzeta})K_{11}(u,\hat{\bzeta})}
{K_{00}(u,\hat{\bzeta})K_{20}(u,\hat{\bzeta})-K_{10}^2(u,\hat
{\bzeta})}.
\]
If the density function $f_{\Lambda_0}$ of $\Lambda_0$ is positive
and the derivative of $E(\ve^2|\Lambda_0=u)$ exists, then we can
further demonstrate that $(nh)^{1/2}\{\wideh\eta(u)-\eta(u)-1/\break2
k_2\eta^{\prime\prime}(u) h^2\}$ converges to a normal distribution $N(0,
\sigma_{\eta,u}^2)$, where $\sigma_{\eta,u}^2=f^{-1}_{\Lambda
_0}(u)\int
K^2(t)\,dt\times  E(\ve^2|\Lambda_0=u)$, $k_2=\int K(t)t^2\,dt$ and $\eta^{\prime\prime}$ is
the second derivative of $\eta$. It is also noteworthy that one
usually either introduces a trimming function or adds a ridge parameter
[see Seifert and Gassert (\citeyear{seigas1996})] when the denominator of $\wideh\eta$ closes
to zero.

\begin{remark}\label{re:1}
Condition (v) indicates that Theorem~\ref{th:th1} is applicable for a
reasonable range of bandwidths. Numerical studies confirm it; our
results remain stable by employing various bandwidths around the
optimal bandwidth selected by cross-validation, in particular when the
sample size becomes large.
\end{remark}

\subsection{Penalized profile least-squares estimator}\label{sec:pseu}

In practice, the true model is often unknown a priori. An underfitted
model can yield biased estimates and predicted values, while an
overfitted model can degrade the efficiency of the parameter estimates
and predictions. This motivates us to apply the penalized least-squares
approach to simultaneously estimate parameters and select important
variables. To this end, we consider a penalized profile least-squares
function
%
\begin{equation}\label{eq:2.5}
{\mathcal L}_P(\bzeta)= \frac{1}{2} Q(\bzeta)+n\sum_{j=1}^{p}
p_{\lambda_{1j}}(|\alpha_j|) +n\sum_{k=1}^{q}
p_{\lambda_{2k}}(|\beta_k|),
\end{equation}
where $p_{\lambda}(\cdot)$ is a penalty function with a regularization
parameter $\lambda$. Throughout this paper, we allow different elements
of $\balpha$ and $\bbeta$ to have different penalty functions with
different regularization parameters. For the purpose of selecting
$X$-variables only, we simply set $p_{\lambda_{1j}}(\cdot)=0$ and the
resulting penalized profile least-squares function becomes
%
\begin{equation}\label{eq:2.6}
{\mathcal L}_P(\bzeta)= \frac{1}{2} Q(\bzeta)+n\sum_{k=1}^{q}
p_{\lambda_{2k}}(|\beta_j|).\vadjust{\goodbreak}
\end{equation}
Similarly, if we are
only interested in selecting $Z$-variables, then we set
$p_{\lambda_{2k}}(\cdot)=0$ so that
\begin{equation}\label{eq:2.7}
{\mathcal L}_P(\bzeta)= \frac{1}{2}
Q(\bzeta)+n\sum_{j=1}^{p} p_{\lambda_{1j}}(|\alpha_j|).
\end{equation}

There are various penalty functions available in the literature. To
obtain the oracle property of Fan and Li (\citeyear{fanli2001}), we adopt their SCAD
penalty, whose first derivative is
\[
p_{\lambda}'(\theta)= \lambda
\biggl\{I(\theta\le\lambda)+\frac{(a\lambda-\theta)_+}{(a-1)\lambda
}I(\theta
>\lambda)\biggr\},
\]
and where $p_{\lambda}(0)=0$, $a=3.7$ and $(t)_{+}=tI\{t>0\}$ is the
hinge loss function. For the given tuning parameters, we obtain the
penalized estimators by minimizing ${\mathcal L}_P(\balpha,\bbeta)$ with
respect to $\balpha$ and $\bbeta$. For the sake of simplicity, we
denote the resulting estimators by $\wideh{\balpha}_{\lambda_1}$ and
$\wideh{\bbeta}_{\lambda_2}$.

In what follows, we study the theoretical properties of the penalized
profile least-squares estimators with the SCAD penalty. Without loss of
generality, it is assumed that the correct model has regression
coefficients $\balpha_0=(\balpha_{10}^{\T},\balpha_{20}^{\T})^{\T}$ and
$\bbeta_0=(\bbeta_{10}^\T,\bbeta_{20}^{\T})^{\T}$, where $\balpha_{10}$ and
$\bbeta_{10}$ are $p_0\times1$ and $q_0\times1$ nonzero components
of $\balpha_0$ and $\bbeta_0$, respectively, and $\balpha_{20}$ and
$\bbeta_{20}$ are $(p-p_0)\times1$ and $(q- q_0)\times1$ vectors with
zeros. In addition, we define $Z_1$ and $X_1$ in such a way that they
consist of the first $p_0$ and $q_0$ elements of $Z$ and $X$,
respectively. We define $\wideti Z_1$ and $\wideti X_1$ analogously. Finally,
we use the following notation for \vspace*{1pt}simplicity: $\GXX=E(\wideti
X_1^{\otimes2})$, $\GXZ=E\{\wideti Z_1 \wideti X_1^\T\eta'(\Lambda)\}$ and
$\GZZ=E\{\wideti Z_1\eta'(\Lambda)\}^{\otimes2}$.\vspace*{-2pt}

\begin{theorem}\label{th:rootn}
Under the regularity conditions \textup{(i)--(vi)}, if, for all $k$ and $j$,
$\lambda_{1j}\to0$, $\sqrt{n}\lambda_{1j}\to\infty$, $\lambda
_{2k}\to
0$ and $\sqrt{n}\lambda_{2k}\to\infty$ as $n\to\infty$, with
probability tending to one, then the penalized estimators
$\hat{\balpha}_{\lambda_1}=(\hat{\balpha}_{1{\lambda_1}}^\T,\hat
{\balpha}_{2{\lambda_1}}^\T)^\T$
and
$\hat{\bbeta}_{\lambda_2}=(\hat{\bbeta}{}^\T_{1{\lambda_2}},\hat
{\bbeta}{}^\T_{2{\lambda_2}})^\T$
satisfy:
\begin{longlist}
\item[(a)] $\hat{\balpha}_{2{\lambda_1}}=0$ and
$\hat{\bbeta}_{2{\lambda_2}}=0$;
\item[(b)] $\sqrt n
(\hat{\balpha}_{1{\lambda_1}}-\balpha_{10})\to N\{0,
\sigma^2(\GZZ-\GXZ\GXX^{-1}\GXZ^\T)^{-1}\}$ and\break $\sqrt n
(\hat{\bbeta}_{1{\lambda_2}}-\bbeta_{10})\to
N\{0,\sigma^2(\GXX-\GXZ^\T\GZZ^{-1}\GXZ)^{-1}\}$.\vspace*{-2pt}
\end{longlist}
\end{theorem}

Analogous results can be established for those parameter estimators
obtained via the penalized functions (\ref{eq:2.6}) or (\ref{eq:2.7}).

Kong and Xia (\citeyear{konxia2007}) developed a cross-validation-based variable
selection procedure in single-index models and observed the difference
between the popular leave-$m$-out cross-validation method in the variable
selection of the linear model and the single-index model. It is
noteworthy that we need two tuning parameters in (\ref{eq:2.5}),
$\lambda_1$ and $\lambda_2$, imposed on the linear part and the
single-index part,\vadjust{\goodbreak} respectively. They can be in the same order in the
large-sample sense, on the basis of the assumptions in Theorem
\ref{th:rootn}, although the distinction between these two tuning
parameters in practice is anticipated, as will be discussed below.
Theorem~\ref{th:rootn} indicates that the proposed variable selection
procedure possesses the oracle property. However, this attractive
feature relies on the tuning parameters. To this end, we adopt Wang, Li
and Tsai's (\citeyear{wanlitsai2007}) BIC selector to choose the regularization parameters
$\lambda_{1j}$ and $\lambda_{2k}$. Because it is computationally
expensive to minimize BIC, defined below, with respect to the
$(p+q)$-dimensional regularization parameters, we follow the approach
of Fan and Li (\citeyear{fanli2004}) to set $\lambda_{1j}=\lambda
\SE(\hat{\alpha}_j^u)$ and $\lambda_{2k}=\lambda\SE(\hat{\beta}_k^u)$,
where $\lambda$ is the tuning parameter, and $\SE(\hat{\alpha}_j^u)$
and $\SE(\hat{\beta}_k^u)$ are the standard errors of the unpenalized
profile least-squares estimators of $\alpha_j$ and $\beta_k$,
respectively, for $j=1,\ldots, p$ and $k=1,\ldots, q$. Let the
resulting SCAD estimators be $\hat{\balpha}_{\lambda}$ and
$\hat{\bbeta}_{\lambda}$. We then select $\lambda$ by minimizing the
objective function\looseness=-1
%
\begin{equation} \label{criterion:BIC}
\BIC(\lambda)
= \log\{(\mse(\lambda)\} +\bigl(\log(n)/n\bigr)\df_{\lambda},
\end{equation}\looseness=0
where $\mse(\lambda)= n^{-1} \sum_{i=1}^n
\{Y_i-\hat{\eta}(Z_i^{\T}\hat{\balpha}_{\lambda})-X_i^{\T}\hat
{\bbeta}_{\lambda}\}^2$
and $\df_{\lambda}$ is the number of nonzero coefficients of both
$\hat{\balpha}_{\lambda}$ and $\hat{\bbeta}_{\lambda}$. More specifically,
we choose $\lambda$ to be the minimizer among a set of grid points over
bounded interval $[0,\lambda_{\mathrm{max}}]$, where $\lambda_{\mathrm{max}}/\sqrt
{n}\to
0$ as $n\to\infty$. The resulting optimal tuning parameter is denoted
by~$\hat\lambda$. In practice, a plot of the BIC($\lambda$) against
$\lambda$ can be used to determine an appropriate $\lambda_{\mathrm{max}}$ to
ensure that the BIC$(\lambda)$ reaches its minimum around the middle of
the range of $\lambda$. The grid points for $\lambda$ are then taken to
be evenly distributed over $[0,\lambda_{\mathrm{max}}]$ so that they are chosen
to be fine enough to avoid multiple minimizers of BIC($\lambda$). In
our numerical studies, the range for $\lambda$ and the number of grid
points are set in the same manner as those in Wang, Li and Tsai (\citeyear{wanlitsai2007})
and Zhang, Li and Tsai (\citeyear{zhatsai2010}). Based on our limited experience, the
resulting estimate of $\lambda$ is quite stable with respect to the
number of grid points when they are sufficiently fine.

To investigate the theoretical properties of the BIC selector, we
denote by $S=\{j_1,\ldots,j_{d}\}$ the set of the indices of the
covariates in the given candidate model, which contains indices of both
$X$ and $Z$. In addition, let $S_{T}$ be the true model,
$S_{F}$ be the full model and $S_{\lambda}$ be the set of the indices
of the covariates selected by the SCAD procedure with tuning parameter
$\lambda$. For a given candidate model $S$ with parameter vectors
$\balpha_s$ and $\bbeta_{s}$, let $\hat{\balpha}_s$ and
$\hat{\bbeta}_s$ be the corresponding profile least-squares
estimators. Then, define $\sigma_n^2(S)=n^{-1} \sum_{i=1}^n
\{Y_i-\hat{\eta}(Z_i^{\T}\hat{\balpha}_s) -X_i^{\T}\hbbeta_s\}
^2$ and
further assume that:
\begin{enumerate}
\item[(A)] for any $S\subset S_F$, $\sigma_n^2(S)\to\sigma^2(S)$ in
probability for some $\sigma^2(S)>0$;
\item[(B)] for any $S\not\supset
S_T$, we have $\sigma^2(S)>\sigma^2(S_T)$.
\end{enumerate}
It is noteworthy that (A) and (B) are the standard conditions for
investigating parameter estimation under model misspecification [e.g.,
see Wang, Li and Tsai~(\citeyear{wanlitsai2007})].\vadjust{\goodbreak}

We next present the asymptotic property of the BIC-type tuning
parameter selector.\vspace*{-2pt}
\begin{theorem}\label{th:BIC}
Under conditions \textup{(A)}, \textup{(B)} and the regularity conditions \textup{(i)--(vi)}, we
have
%
\begin{equation}\label{eq:2.9}
P(S_{\widehat{\lambda}}=S_{T})\to1.\vspace*{-2pt}
\end{equation}
\end{theorem}

Theorem~\ref{th:BIC} demonstrates that the BIC tuning
parameter selector enables us to select the true model consistently.\vspace*{-2pt}

\section{Hypothesis tests}\label{sec:test}

Applying the estimation method described in the previous section, we
propose two hypothesis tests. The first is for general hypothesis
testing of regression parameters and the second is for testing the
nonparametric function.\vspace*{-2pt}

\subsection{Testing parametric components}\label{subsec:test1}

Consider the general linear hypothesis
%
\begin{equation}\label{eq:3.1}
H_0\dvtx \bA\bzeta=\bdelta\quad\mbox{versus}\quad H_1\dvtx \bA\bzeta\neq\bdelta,
\end{equation}
where
$\bA$ is a known $m\times(p+q)$ full-rank matrix and $\bdelta$ is an
$m\times1$ vector. A~simple example of (\ref{eq:3.1}) is to test
whether some elements of $\balpha$ and $\bbeta$ are zero; that is,
\[
H_0\dvtx \alpha_{i_1}=\cdots=\alpha_{i_k}=0\mbox{ and }\beta_{j_1}=\cdots=\beta_{j_l}=0
\]
versus
\[
H_1\dvtx \mbox{not all }\alpha_{i_1},\ldots,\alpha_{i_k}
\mbox{ and }\beta_{j_1}, \ldots,\beta_{j_l}\mbox{ are equal to }0.
\]

Under $H_0$ and $H_1$, let $\bzeta_0=(\balpha_0^\T,\bbeta_0^\T)^\T
$ and
$\bzeta_1=(\balpha_1^\T,\bbeta_1^\T)^\T$ be the corresponding parameter
vectors, and let $\Omega_0$ and $\Omega_1$ be the parameter spaces of
$\bzeta_0$ and $\bzeta_1$, respectively. It is noteworthy that this is
a slight abuse of notation because $\bzeta_0$ was previously used to
denote the true value of $\bzeta$ in Section~\ref{sec:plse}. Furthermore, define
\[
Q(H_0)=\inf_{\Omega_0}
Q(\bzeta_0)=\sum^n_{i=1}\{Y_i-{\bX}_i^\T\hat{\bbeta}_0-\hat{\eta
}(Z_i^{\T}\hat{\balpha}_0,
\hat\bzeta_0)\}^2
\]
and
\[
Q(H_1)=\inf_{\Omega_1}Q(\bzeta_1)=\sum^n_{i=1}\{Y_i-{\bX}_i^\T
\hat{\bbeta}_1-\hat{\eta}(Z_i^{\T}\hat\balpha_1,
\hat\bzeta_1)\}^2,
\]
where $\{\hat{\bbeta}_0, \hat\bzeta_0\}$ and
$\{\hat{\bbeta}_1, \hat\bzeta_1\}$ are the profile least-squares
estimators of $\{\bbeta_0, \bzeta_0\}$ and $\{\bbeta_1, \bzeta_1\}$,
respectively, and $\hat{\eta}$ is the nonparametric estimator of
$\eta$
obtained via (\ref{eq:2.3}). Subsequently, we propose a test statistic,
\[
T_1 = \frac{n\{Q(H_0)-Q(H_1)\}}{Q(H_1)},
\]
and give its theoretical property below.\vadjust{\goodbreak}

\begin{theorem} \label{th:test1}
Assume that the regularity conditions \textup{(i)}--\textup{(vi)} hold. Then:
\begin{longlist}
\item[(a)] under $H_0$ in (\ref{eq:3.1}), $T_1 \rightarrow\chi_m^2;$
\item[(b)] under $H_1$ in (\ref{eq:3.1}), $T_1$ converges to a
noncentral chi-squared distribution with $m$ degrees of freedom and
noncentrality parameter $\phi=\break\lim_{n\to\infty}
n\sigma^{-2}(\bA\bzeta-\bdelta)^\T(\bA\bfD^{-1}\bA^\T)^{-1}(\bA\bzeta-\bdelta)$,
where $\bfD$ is defined as in Theorem~\ref{th:th1}.
\end{longlist}
\end{theorem}

Analogously, \vspace*{1pt}we are able to construct the Wald test,
$W_n=(\bA\hat{\bzeta}-\break\bdelta)^\T(\bA\hat{\bfD}^{-1}\times\bA^\T
)^{-1}(\bA\hat{\bzeta}-\bdelta)$, and demonstrate that $W_n$ and $T_1$
have the same asymptotic distribution.

\subsection{Testing the nonparametric component}\label{subsec:test2}

The nonparametric estimate of $\eta(\cdot)$ provides us with
descriptive and graphical information for exploratory data analysis.
Using this information, it is possible to formulate a parametric model
that takes into account the features that emerged from the preliminary
analysis. To this end, we introduce a goodness-of-fit test to assess
the appropriateness of a proposed parametric model. Without loss of
generality, we consider a simple linear model under the null
hypothesis. Accordingly, the null and alternative hypotheses are given
as follows:
\begin{equation}\label{eq:3.2}
\hspace*{20pt}H_0\dvtx\eta(u) = \theta_0+\theta_1u \quad\mbox{versus}\quad
H_1\dvtx\eta(u)\ne\theta_0 + \theta_1u\qquad\mbox{for some }
u,\nonumber
\end{equation}
where $\theta_0$ and $\theta_1$ are unknown constant parameters.

Under $H_1$, let $\wideh\balpha$, $\hbbeta$ and $\wideh{\eta}$ be the
corresponding profile least-squares and nonparametric estimators of
$\balpha$, $\bbeta$ and $\eta$, respectively. Under $H_0$, we use the
same parametric estimators $\wideh\balpha$ and $\hbbeta$ as those
obtained under $H_1$, while the estimator of $\eta$ is
$\tilde{\eta}(u)=\hat\theta_0 + \hat\theta_1 u$, where $\hat
\theta_0$
and $\hat\theta_1$ are the ordinary least-squares estimators of
$\theta_0$ and $\theta_1$, respectively, by fitting
$Y_i-X_i^\T\wideh\bbeta$ versus $Z_i^\T\wideh\balpha$. The resulting residual
sums of squares under the null and alternative hypotheses are then
\[
\rss(H_0)=
\sum_{i=1}^n\{Y_i-\tilde{\eta}(Z_i^\T\wideh\balpha)-X_i^\T\wideh\bbeta\}^2
\]
 and
\[
\rss(H_1)= \sum_{i=1}^n
\{Y_i-\wideh\eta(Z_i^\T\wideh\balpha)-X_i^\T\wideh\bbeta\}^2.
\]

To test the null hypothesis, we consider the following generalized
$F$-test:
\[
T_2=\frac{r_K}2\frac{n\{\rss(H_0)-\rss(H_1)\}}{\rss(H_1)},
\]
where $r_K=\{K(0)-0.5\int K^2(u)\,du\}\{\int\{K(u)-0.5K*K(u)\}\,du\}^{-1}$
and $K*K$ denotes the convolution of $K$. The theoretical property of
$T_2$ is given below.

\begin{theorem}\label{th:glrt}
Assume that the regularity conditions \textup{(i)}--\textup{(vi)} hold. Under $H_0$ in
(\ref{eq:3.2}), $T_2$ has an asymptotic $\chi^2$ distribution with
$df_n$ degrees of freedom,  where $df_n=r_K |{\mathcal
U}|\{K(0)-0.5\int K^2(u)\,du\}/h$, $|{\mathcal U}|$ stands for the length of
${\mathcal U}$ and $\mathcal U$ is defined in regularity condition \textup{(i)}.
\end{theorem}

The above theorem unveils the Wilks phenomenon for the
partially linear single-index model. Furthermore, we can obtain an
analogous result when the simple linear model under $H_0$ is replaced
by a multiple regression model.

\section{Simulation studies}\label{sec:nume}

In this section, we present four Monte Carlo studies which evaluate the
finite-sample performance of the proposed estimation and testing
methods. The first two examples illustrate the performance of the
profile least-squares estimator and the SCAD-based variable selection
procedure proposed in Sections~\ref{sec:plse} and~\ref{sec:pseu},
respectively. The next two examples explore the performance of the test
statistics developed in Sections~\ref{subsec:test1} and~\ref{subsec:test2}.

\begin{example}\label{ex:1}
We generated $500$ realizations, each consisting of $n=50, 100$
and $200$ observations, from each of the following models:
%
\begin{eqnarray}\label{eq:4.1}
y&=&4\bigl\{(z_1+z_2-1)/{\sqrt 2}\bigr\}^2+4+0.2\ve; \\\label{eq:4.2}
y&=&\sin\bigl[\bigl\{(z_1+z_2+z_3)/{\sqrt 3}-a\bigr\}\pi/(b-a)\bigr]+\beta X+0.1\ve,
\end{eqnarray}
where $z_1, z_2, z_3$ are independent and
uniformly distributed on $[0,1]$, $X=0$ for the odd numbered
observations and $X=1$ for the even numbered observations, $\ve$~has
the standard normal distribution, $a=0.3912$ and $b=1.3409$. The
resulting parameters of models (\ref{eq:4.1}) and (\ref{eq:4.2}) are
$\balpha=(\alpha_1,\alpha_2)^{\T}=(0.7071,0.7071)^{\T}$ and
$(\balpha,\beta)=(\alpha_1,\alpha_2,\alpha_3,\beta)^{\T}=(0.5774,0.5774,
0.5774, 0.3)^{\T}$ with $\phi_0=\arccos(\alpha_1)=0.7854$ as in Xia and
H\"ardle (\citeyear{xiahar2006}) [or $\pi/4$ in H\"ardle, Hall and Ichimura (\citeyear{harhalich1993}), page
165].
\end{example}

Model (\ref{eq:4.1}) was analyzed by H\"ardle, Hall and Ichimura (\citeyear{harhalich1993})
and Xia and H\"ardle (\citeyear{xiahar2006}), while model (\ref{eq:4.2}) was investigated by Carroll
et al. (\citeyear{carfangijwan1997}) and Xia and H\"ardle (\citeyear{xiahar2006}).
For both models, Xia and
H\"ardle claimed that their MAVE approach outperforms those of
H\"ardle, Hall and Ichimura (\citeyear{harhalich1993}) and Carroll et al. (\citeyear{carfangijwan1997}),
respectively. It is
therefore of interest to compare the profile least-squares (abbreviated
to PrLS) method with MAVE. Tables~\ref{ta:1} and~\ref{ta:2} present the
results for models (\ref{eq:4.1}) and  (\ref{eq:4.2}), respectively.
Both tables indicate that the profile least-squares method  yields
accurate estimates and that the mean squared error becomes smaller as
the sample size gets larger, which is consistent with the theoretical
finding. Furthermore, Table~\ref{ta:1} shows that the mean squared
errors of $\wideh\balpha$ and $\wideh\bbeta$ are smaller than those
computed via the  MAVE method [see Table 1 of Xia and H\"ardle (\citeyear{xiahar2006})].
Table~\ref{ta:2} suggests that the biases\vadjust{\goodbreak} and their associated mean
squared errors of $\wideh\balpha$ and $\wideh\bbeta$ are comparable to
those calculated via MAVE. In summary, the Monte Carlo studies indicate
that PrLS performs well.  Because the penalized MAVE is not easy to
obtain, we study only the penalized PrLS estimates in the following
example.

\begin{table}
\tablewidth=310pt
\caption{Simulation results for Example \protect\ref{ex:1}: the profile
least-squares estimates (PrLS) and their corresponding mean squared
errors ($\times 10^{-4}$) for model (\protect\ref{eq:4.1})}\label{ta:1}
\begin{tabular*}{310pt}{@{\extracolsep{\fill}}lccc@{}}
\hline
$\bolds{n}$ & $\bolds{\alpha_{1}\,(=0.7071)}$  & $\bolds{\alpha_{2}\,(=0.7071)}$ & $\bolds{\phi_0\,(=0.7854)}$ \\
\hline
\phantom{0}50& \phantom{0}0.7053 (21.5274)& \phantom{0}0.7059 (21.3398)& 0.7859 (42.9158)\\
100& 0.7054 (8.5874)& 0.7076 (8.4175)& 0.7869 (17.0116)\\
200& 0.7067 (4.4636)& 0.7069 (4.4287)& 0.7856 (8.8942)\phantom{0}\\
\hline
\end{tabular*}
\end{table}

\begin{table}[b]
\tablewidth=315pt
\caption{Simulation results for Example \protect\ref{ex:1}: the profile
least-squares estimates (PrLS) and their corresponding mean squared
errors ($\times 10^{-4}$) for model (\protect\ref{eq:4.2})}\label{ta:2}
\begin{tabular*}{315pt}{@{\extracolsep{\fill}}lcccc@{}}
\hline
$\bolds{n}$ & $\bolds{\alpha_{1}\,(=0.5774)}$  & $\bolds{\alpha_{2}\,(=0.5774)}$ & $\bolds{\alpha_{3}\,(=0.5774)}$ & $\bolds{\beta\,(=0.3)}$ \\
\hline
\phantom{0}50& 0.5753 (5.8336)& 0.5771 (5.5685)& 0.5782 (5.9245)& \phantom{0}0.2923 (11.4582)\\
100& 0.5776 (2.5009)& 0.5774 (2.4606)& 0.5764 (2.4770)& 0.3000 (4.7030)\\
200& 0.5782 (1.1533)& 0.5771 (1.0852)& 0.5764 (1.2483)& 0.3004 (2.2026)\\
\hline
\end{tabular*}
\end{table}


\begin{example}\label{ex:2}
We simulated $500$ realizations, each consisting of $n=100$ and
$200$ random samples,  from model\vspace*{1pt} (\ref{eq:4.2}) with $\sigma=0.1$ and
$0.25$, respectively. The mean function has coefficients
$\balpha=(1,3,1.5,0.5,0,0,0,0)^{\T}/\sqrt{12.5}$ and
$\bbeta=(3,2,0,0,0,1.5,0,0.2,0.3,0.15,0,0)^{\T}$. To assess the
robustness of estimates, we further generate the linear and nonlinear
covariates from the following three scenarios: (i) the covariate
vectors $\bfX$ and $\bfZ$ have $12$ and $8$ elements, respectively,
which are independent and  uniformly distributed on $[0,1]$;
(ii) the covariate vector $\bfX$ has $12$ elements; the first 5 and
last 5 elements are independent and standard normally distributed,
while the 6th and 7th elements are independently Bernoulli distributed
with success probability $0.5$; the covariate vector $\bfZ$ has $8$
elements, which are independent and  standard normally distributed;
(iii) a covariate vector $W$ was generated from a 12-dimensional normal
distribution with mean 0 and variance 0.25; the correlation between
$w_i$ and $w_j$ is $\rho^{|i-j|}$ with $\rho=0.4$;  then the covariate
vector $\bfX$ =
$W+\{1.5\exp{(1.5z_1)},5z_1,5\sqrt{z_2},3z_1+z_2^2,0,0,0,0,0,0,0,0\}^{\T}$;
moreover, the covariate vector $\bfZ$ has $8$ elements, which are
independent and  uniformly distributed on $[0,1]$.
\end{example}

\begin{table}
\tabcolsep=0pt
\caption{Simulation results for Example \protect\ref{ex:2}. S-AIC:
SCAD(AIC); S-BIC: SCAD(BIC). MRME:~median of relative model error; C:
the average number of the true zero coefficients that were
correctly set to zero; I: the average number of the truly nonzero~coefficients that were incorrectly set to zero}\label{ta:i}
\begin{tabular*}{\textwidth}{@{\extracolsep{\fill}}lccccccccccccc@{}}
\hline
&&\multicolumn{6}{c}{$\bolds{\sigma=0.1}$}&\multicolumn{6}{c@{}}{$\bolds{\sigma=0.25}$}\\ [-7pt]
&&\multicolumn{6}{c}{\hrulefill}&\multicolumn{6}{c@{}}{\hrulefill}\\
&&\multicolumn{3}{c}{$\balpha$}&\multicolumn{3}{c}{$\bbeta$} &  \multicolumn{3}{c}{$\balpha$}  &\multicolumn{3}{c@{}}{$\bbeta$}
\\ [-7pt]
&&\multicolumn{3}{c}{\hrulefill}&\multicolumn{3}{c}{\hrulefill}&\multicolumn{3}{c}{\hrulefill}&\multicolumn{3}{c@{}}{\hrulefill}\\
 $\bolds{n}$ & &  \textbf{MRME} &  \textbf{C}  & \textbf{I}  &  \textbf{MRME} &  \textbf{C}
 & \textbf{I} &  \textbf{MRME} &  \textbf{C}  & \textbf{I} &  \textbf{MRME} &  \textbf{C}  & \textbf{I}\\
\hline
&& \multicolumn{12}{c@{}}{scenario (i)} \\
100& Oracle& 0.26& 4\phantom{00,}& 0\phantom{00,}& 0.28& 6\phantom{00,}& 0\phantom{00,}& 0.2\phantom{0}& 4\phantom{00,}& 0\phantom{00,}& 0.27& 6\phantom{00,}& 0\phantom{00,}\\
& S-BIC& 0.37& 3.60& 0.08& 0.91& 5.32& 0.29& 0.73& 3.29& 0.30& 0.86& 4.91& 1.02\\
& S-AIC& 0.66& 3.08& 0.05& 0.97& 4.12& 0.15& 0.75& 2.70& 0.11& 0.91& 4.02&
0.68\\ [3pt]
200& Oracle& 0.27& 4\phantom{00,}& 0\phantom{00,}& 0.34& 6\phantom{00,}& 0\phantom{00,}& 0.31& 4\phantom{00,}& 0\phantom{00,}& 0.39& 6\phantom{00,}& 0\phantom{00,}\\
& S-BIC& 0.33& 3.89& 0.02& 0.85& 5.55& 0.02& 0.36& 3.86& 0.03& 0.94& 5.50& 0.57\\
& S-AIC& 0.60& 3.39& 0.01& 0.92& 4.49& 0.01& 0.62& 3.29& 0.01& 0.93& 4.43&
0.23\\ [6pt]
&& \multicolumn{12}{c@{}}{scenario (ii)}\\
100& Oracle& 0.29& 4\phantom{00,}& 0\phantom{00,}& 0.35& 6\phantom{00,}& 0\phantom{00,}& 0.24& 4\phantom{00,}& 0\phantom{00,}& 0.24& 6\phantom{00,}& 0\phantom{00,}\\
& S-BIC& 0.36& 3.75& 0.05& 0.88& 5.44& 0.19& 0.66& 3.47& 0.27& 0.94& 5.11& 1.07\\
& S-AIC& 0.65& 3.26& 0.02& 0.91& 4.35& 0.07& 0.70& 2.86& 0.09& 0.96& 4.04&
0.67\\ [3pt]
200& Oracle& 0.31& 4\phantom{00,}& 0\phantom{00,}& 0.36& 6\phantom{00,}& 0\phantom{00,}& 0.32& 4\phantom{00,}& 0\phantom{00,}& 0.3& 6\phantom{00,}& 0\phantom{00,}\\
& S-BIC& 0.36& 3.91& 0.01& 0.79& 5.64& 0.01& 0.40& 3.87& 0.03& 0.85& 5.53& 0.50\\
& S-AIC& 0.62& 3.32& 0\phantom{00,}& 0.93& 4.51& 0.01& 0.72& 3.29& 0.01& 0.97& 4.45&
0.19\\ [6pt]
&& \multicolumn{12}{c}{scenario (iii)}\\
100& Oracle& 0.28& 4\phantom{00,}& 0\phantom{00,}& 0.15& 6\phantom{00,}& 0\phantom{00,}& 0.19& 4\phantom{00,}& 0\phantom{00,}& 0.18& 6\phantom{00,}& 0\phantom{00,}\\
& S-BIC& 0.48& 3.67& 0.03& 0.82& 5.24& 0.05& 0.50& 3.35& 0.21& 0.85& 4.99& 0.56\\
& S-AIC& 0.74& 3.09& 0.02& 0.94& 4.35& 0.03& 0.71& 2.70& 0.12& 0.98& 4.17&
0.32\\ [3pt]
200& Oracle& 0.32& 4\phantom{00,}& 0\phantom{00,}& 0.18& 6\phantom{00,}& 0\phantom{00,}& 0.29& 4\phantom{00,}& 0\phantom{00,}& 0.17& 6\phantom{00,}& 0\phantom{00,}\\
& S-BIC& 0.39& 3.89& 0\phantom{00,}& 0.73& 5.52& 0.01& 0.39& 3.80& 0.04& 0.83& 5.29& 0.12\\
& S-AIC& 0.68& 3.30& 0\phantom{00,}& 0.84& 4.54& 0\phantom{00,}& 0.66& 3.13& 0.02& 0.85& 4.48& 0.03\\
\hline
\end{tabular*}
\end{table}

Based on the above model settings, we next explore the performance of
the penalized profile least-squares approach via SCAD-BIC. Because the
Akaike information criterion [Akaike (\citeyear{aka1973})] has been commonly used for
classical variable selections, we also study SCAD-AIC by replacing
$\log(n)$ in (\ref{criterion:BIC}) with $2(p+q)$. To assess the
performance, we  consider Fan and Li's (\citeyear{fanli2001}) median of relative model
error (MRME), where the relative model error is defined as
$\mathrm{RME}=\mathrm{ME}/\mathrm{ME}_{S_{F}}$,
ME is defined as $E(Z^{\T}\hat\balpha_{\hat\lambda}-Z^{\T}\balpha)^2$
for $\hat\balpha_{\hat\lambda}$ and
$E(X^{\T}\hat\bbeta_{\hat\lambda}-X^{\T}\bbeta)^2$
for $\hat\bbeta_{\hat\lambda}$, and ME$_{S_{F}}$
is the corresponding model error calculated by fitting the data
with the full model via the unpenalized estimates. In
addition, we calculate the average number of the  true zero
coefficients that were correctly set to zero and the average number of
the truly nonzero coefficients that were incorrectly set to zero.
Table~\ref{ta:i} shows that the SCAD-BIC outperforms SCAD-AIC in terms
of model error measures. Moreover, SCAD-BIC has a much better rate of
correctly identifying the true submodel than that of SCAD-AIC, although
it sometimes shrinks  small  nonzero coefficients to zero.
Unsurprisingly, SCAD-BIC improves as the signal gets stronger and the
sample size becomes larger, which corroborates our theoretical
findings.

\begin{figure}

\includegraphics{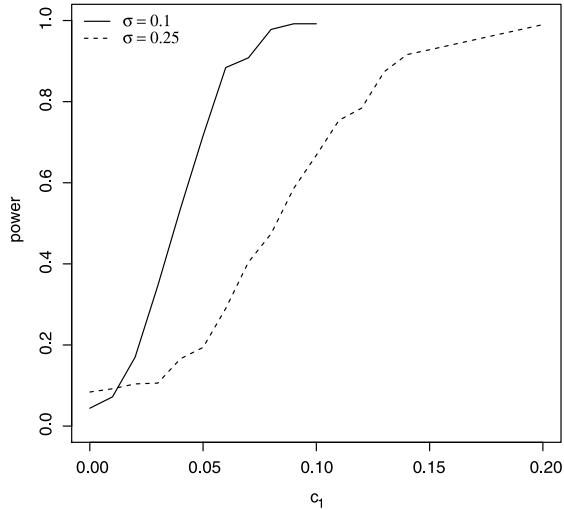}

\caption{Power function of the test statistic
$T_1$.} \label{fig:parpow}
\end{figure}

\begin{example}\label{ex:3}
To study the finite-sample performance of the test statistic $T_1$
in Section~\ref{subsec:test1}, we consider the same model as that of
scenario (i) in Example~\ref{ex:2}. Due to the model's parameter
setting, we naturally consider the following null and alternative
hypotheses:
\[
H_0\dvtx \beta_{3}=\beta_{4}=\beta_{5}=\beta_{7}=0 \quad\mbox{versus}\quad
H_1\dvtx
\beta_{3}=\beta_{4}=\beta_{5}=\beta_{7}=c_1,
\]
where $c_1$ ranges from 0 to 0.1 with increment 0.01  for $\sigma=0.1$,
whereas $c_1$ is the value from a set of $\{0, 0.01, 0.02, \ldots,
0.09, 0.15, 0.2\}$ for $\sigma=0.25$. In addition, $500$ realizations
were generated with $n=200$ to calculate the size and power of $T_1$.
Figure~\ref{fig:parpow} depicts the power function versus $c_1$. It
shows that the empirical size at $c_1=0$ is very close to the nominal
level 0.05. Furthermore,  the power of the test is greater than $0.95$
as $c_1$ increases to $0.05$  and $0.15$, respectively, when
$\sigma=0.1$ and $\sigma=0.25$. It is not surprising that the power
increases as the signal gets stronger. In summary, $T_1$ not only
controls the size well, but is also a powerful test.
\end{example}

\begin{figure}

\includegraphics{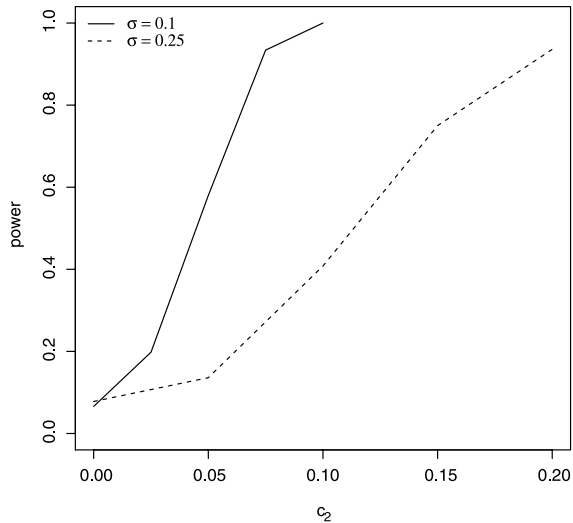}

\caption{\small Power function of the test statistic
$T_2$.} \label{fig:nonparpow}
\end{figure}

\begin{example}\label{ex:4}
To examine the performance of the test statistic $T_2$ in Section~\ref{subsec:test2},
we generated $500$ realizations from the model
given below with $n=200$.
%
\begin{equation}\label{eq:4.3}
y=\eta\bigl\{(z_1+z_2+z_3)/{\sqrt 3}\bigr\}-
0.5x_1 + 0.3x_2+\sigma\ve,
\end{equation}
where $\sigma=0.1$ and
$\sigma=0.25$, respectively. We then consider the following hypotheses:
\[
H_0\dvtx\eta(u)= u \quad\mbox{versus}\quad H_1\dvtx\eta(u)=
c_2\sin\{\pi(u-a)/(b-a)\}+u,
\]
where $c_2$ ranges from 0 to 0.1 with
increment 0.025 for $\sigma=0.1$, while $c_2$ ranges from 0 to 0.2 with
increment 0.05 for $\sigma=0.25$. Figure~\ref{fig:nonparpow}
demonstrates that the empirical size at $c_2=0$ is very close to the
nominal level 0.05. Furthermore,  the power of the test is greater than
$0.95$  as $c_2$ increases to $0.075$  and $0.2$, respectively. As
expected, the power increases when the signal becomes stronger.
Although $T_2$ is slightly less powerful than $T_1$, it controls the
size well and is a  reliable test.
\end{example}

\section{Discussion}\label{sec:disc}

In partially linear  single-index models, we propose using the SCAD
approach to  shrink parameters contained in both parametric and
nonpara-\ metric components. The resulting estimators enjoy the oracle
property when the regularization parameters satisfy the proper
conditions. To further exploit SCAD, one could extend the current
results to partially linear multiple-index models\vspace*{1pt} by allowing $\ve$ to
be dependent on $(Z^{\T},X^{\T})^{\T}$. In addition, one could obtain
the SCAD estimator for generalized partially linear single-index
models. Finally, an investigation of partially linear single-index
model selection with  error-prone covariates could also be of interest.
We believe that these efforts would enhance the usefulness of SCAD in
data analysis.

\begin{appendix}
\section*{Appendix: Proofs of theorems}\label{sec:appe}

\subsection{\texorpdfstring{Proof of Theorem \protect\ref{th:th1}}{Proof of Theorem 1}}

Under the conditions of Theorem~\ref{th:th1}, we follow similar arguments to those
used by Ichimura (\citeyear{ich1993}) and show that
$\wideh\bzeta=(\wideh\balpha{}^{\T},\wideh{\bbeta}{}^\T)^{\T}$ is a root-$n$
consistent estimator of $\bzeta$. Because the proof is straightforward,
we do not present it here. We next demonstrate the  asymptotic
normality of $\wideh\bzeta$ by using a general result of Newey (\citeyear{new1994}).

Let $m_x(\Lambda) = E(X|\Lambda)$, $m_z(\Lambda) = E(Z|\Lambda)$ and
$\kappa =\eta'(\Lambda)\{Z-m_z(\Lambda)\}$. In addition, let
%
\begin{equation}\label{eq:psi}
\Psi(m_x,\eta,\kappa,\balpha,\bbeta,Y,Z,X)=(Y-\eta-X^\T\bbeta)
\left\{\matrix{\kappa\vspace*{4pt}\cr
X-m_x(\Lambda)}\right\}.
\end{equation}
For any given $m_x^*$, $\eta^*$ and $\kappa^*$,
define
\begin{eqnarray*}
&& D(m_x^*-m_x,\eta^*-\eta,\kappa^*-\kappa,\balpha, \bbeta, Y,Z,X) \\
&&\qquad =\frac{\partial\Psi}{\partial m_x}(m_x^*-m_x)
+\frac{\partial\Psi}{\partial\eta}(\eta^*-\eta)
+\frac{\partial\Psi}{\partial\kappa}(\kappa^*-\kappa),
\end{eqnarray*}
where the
partial derivatives are the Frechet partial derivatives. After
algebraic simplification, we have
\begin{eqnarray*}
\frac{\partial\Psi}{\partial
m_x}&=&(Y-\eta-X^\T\bbeta) \pmatrix{0\cr-1},\\
\frac{\partial\Psi}{\partial \eta}&=& -\left\{\matrix{\kappa\vspace*{2pt} \cr
X-m_x(\Lambda)}
\right\},\\
\frac{\partial\Psi}{\partial \kappa}&=&(Y-\eta-X^\T\bbeta)
\pmatrix{1\cr 0},
\end{eqnarray*}
where the partial derivatives are zero. Accordingly,
%
\begin{eqnarray}\label{eq:amyeq4}
&& \Vert\Psi(m^*_x,\eta^*,\kappa^*,\balpha, \bbeta,Y,Z,X) -
\Psi(m_x,\eta,\kappa,\balpha, \bbeta,Y,Z,X)\nonumber\\
&& \qquad \quad\hspace*{24pt}{}-D(m_x^*-m_x,\eta^*-\eta,\kappa^*-\kappa,\balpha,
\bbeta,Y,Z,X)\Vert\\
&& \qquad= O( \Vert m_x^* - m_x \Vert^2 + \Vert \eta^* -
\eta \Vert^2 + \Vert \kappa^* - \kappa \Vert^2),\nonumber
\end{eqnarray}
where $\Vert\cdot\Vert$ denotes the Sobolev norm, that is, the supremum norm of the
function itself, as well as its derivatives. Equation (\ref{eq:amyeq4})
is Newey's Assumption~5.1(i). It is also noteworthy that his Assumption
5.2 holds by the expression of $D(\cdot, \cdot, \cdot, \balpha, \bbeta,
Y, Z, X)$. Moreover, the result
\[
E\{D(m_x^*-m_x,\eta^*-\eta,\kappa^*-\kappa,\balpha, \bbeta,Y,Z,X) \} = 0
\]
leads to Newey's Assumption 5.3.

In addition to  Newey's assumptions mentioned above, we need to verify
one more assumption before employing his result. To this end, we
re-express the solution of
 (\ref{eq:2.2}) as
\[
\wideh \eta(u)=\frac1n \sum_{i=1}^n\frac{\{\wideh s_2(u,h)-\wideh
s_1(u,h)(\Lambda_i-u)\}K_h(\Lambda_i-u)(Y_i-X_i^\T\bbeta)} {\wideh
s_2(u,h)\wideh s_0(u,h)-\wideh s_1^2(u,h)},
\]
where $\Lambda_i$ is the
$i$th row of $\Lambda$ and
\[
\wideh s_r(u,h)=\frac1n\sum_{i=1}^n(\Lambda_i-u)^rK_h(\Lambda_i-u) \qquad\mbox{for } r=0,1,2.
\]
Then, let
\begin{eqnarray*}
\wideh m_x(u)&=&\frac1n \sum_{i=1}^n\frac{\{\wideh
s_2(u,h)-\wideh s_1(u,h)(\Lambda_i-u)\}K_h(\Lambda_i-u)X_i^\T}
{\wideh s_2(u,h)\wideh s_0(u,h)-\wideh s_1^2(u,h)},\\
\wideh m_z(u)&=&\frac1n \sum_{i=1}^n\frac{\{\wideh s_2(u,h)-\wideh
s_1(u,h)(\Lambda_i-u)\}K_h(\Lambda_i-u)Z_i^\T} {\wideh s_2(u,h)\wideh
s_0(u,h)-\wideh s_1^2(u,h)}.
\end{eqnarray*}
Applying similar techniques to those
used in Mack and Silverman  (\citeyear{macsil1982}), we obtain the following equations,
which hold uniformly in $u\in{\mathcal U}$:
\begin{eqnarray}\label{eq:unrate}
\wideh\eta(u)-\eta(u)&=&o_p(n^{-1/4}),\qquad
\wideh\eta'(u)-\eta'(u)=o_p(n^{-1/4}),\nonumber\\ [-8pt]\\ [-8pt]
\wideh m_x(u)-m_x(u)&=&o_p(n^{-1/4})\quad\mbox{and}\quad\wideh
m_z(u)-m_z(u)=o_p(n^{-1/4}).\nonumber
\end{eqnarray}
These results imply that $\wideh{\kappa}- \kappa= o_p(n^{-1/4})$. Thus,
Newey's  Assumption~5.1(ii) holds.

After examining Newey's Assumptions 5.1--5.3, we apply his Lemma 5.1
and find that $\wideh\bzeta$ has the same limit distribution as the
solution to the equation
%
\begin{equation}\label{eq:amyeq5}
0 = \sum_{i=1}^n\Psi(m_x,\eta,\kappa,\balpha,\bbeta,Y_i,Z_i,X_i).
\end{equation}
Furthermore, it is easy to show that the
solution to (\ref{eq:amyeq5}) has the same limit distribution as
described in the statement of Theorem~\ref{th:th1}. Hence, we complete
the proof of asymptotic normality.

Finally, we show the efficiency of $\wideh\bzeta$. Let $p_{\ve}(\ve)$ be
the probability density function of $\ve$ and let $p'_{\ve}(\ve)$ be
its first-order derivative with respect to $\ve$. Then, the score
function of $(\balpha, \bbeta)$ is
\[
S_{(\balpha,\bbeta)}=-\sigma^2
\left\{\matrix{\eta'(\Lambda)Z\vspace*{2pt}\cr X}\right\}
\frac{p'_{\ve}(\ve)}{p_{\ve}(\ve)}.
\]

For any given function $g$ of $(Z,X)$, it can be shown that the
nuisance tangent space ${\mathcal P}$, for the three nuisance parameters,
$g_{(Z,X)}(z,x)$, $p_{\ve}(\ve)$ and $\eta(\Lambda)$, is $\{g(Z,X)\dvtx
E(g)=0, E(\ve g)$ is a function of $(Z,X)$ only$\}$. Furthermore, the
orthogonal component of ${\mathcal P}$ is
\[
{\mathcal P}^{\bot}=\{\ve g(Z,X)\dvtx E(g|\Lambda)=0\}.
\]
Subsequently, we apply the approach of  Bickel et al. (\citeyear{bicklaritwel1993})  and
obtain the following semiparametric efficient score function via
equation (\ref{eq:2.4}):
%
\begin{equation}\label{eq:eff}
S_{\mathrm{eff}}=\ve \left\{
\matrix{\eta'(\Lambda)\wideti Z\vspace*{2pt}\cr
\wideti X}\right\}.
\end{equation}
It can be seen that $S_{\mathrm{eff}}\in {\mathcal P}^{\bot}$.

For any $\ve g\in {\mathcal P}^{\bot}$, we have $E(g|\Lambda)=0$.
Accordingly,
\begin{eqnarray*}
&&E\bigl\{\bigl(S_{(\balpha,\bbeta)}-S_{\mathrm{eff}}\bigr)^\T\ve g(Z,X)\bigr\} \\
&&\qquad=E\left( -\sigma^2\left\{
\matrix{\eta'(\Lambda)Z\vspace*{2pt}\cr
X}\right\}^\T g E\left\{\frac{\ve p'_{\ve}(\ve)}{p_{\ve}(\ve)}\right\}\right.\\
&&\qquad\quad\hspace*{15pt}{}- \left.E(\ve^2)\left[\left\{
\matrix{\eta'(\Lambda)Z\vspace*{2pt}\cr X}\right\}^\T g- E\left\{
\matrix{\eta'(\Lambda)E(Z|\Lambda)\vspace*{2pt}\cr E(X|\Lambda)}\right\}^\T g\right]\right).
\end{eqnarray*}
Because $E\{{\ve p'_{\ve}(\ve)}/{p_{\ve}(\ve)}\}=-1$,
it follows that
\[
E\bigl\{\bigl(S_{(\balpha,\bbeta)}-S_{\mathrm{eff}}\bigr)^\T\ve g(Z,X)\bigr\}
=E\bigl[\{\eta'(\Lambda)E(Z^\T|\Lambda), E(X^\T|\Lambda)\}
E(g|\Lambda)\bigr]=0.
\]
That is, $S_{\mathrm{eff}}$ is the projection of
$S_{(\balpha,\bbeta)}$ onto ${\mathcal P}^{\bot}$ and  the estimator
$\wideh\bzeta$ is therefore efficient [see Bickel et al. (\citeyear{bicklaritwel1993})]. We have
thus completed the proof.\looseness=-1

\subsection{\texorpdfstring{Proof of Theorem \protect\ref{th:rootn}}{Proof of Theorem 2}}

To prove this theorem, we consider the following three steps:  Step I
establishes the order of the minimizer $(\wideh\balpha{}^{\T}_{\lambda_1},
\wideh{\bbeta}{}^{\T}_{\lambda_2})^{\T}$ of ${\mathcal L}_P(\balpha,\bbeta)$; Step
II  shows that $(\wideh\balpha{}^{\T}_{\lambda_1},
\wideh{\bbeta}{}^{\T}_{\lambda_2})^{\T}$  attains sparsity; Step III derives
the asymptotic distribution of the penalized estimators.

\textit{Step I}. Let $\gamma_n = n^{-1/2} + a_n+c_n$,
$\bv_1=(v_{11},\ldots,v_{1p})^{\T}$,
$\bv_2=(v_{21},\ldots,v_{2q})^{\T}$ and $\|\bv_1 \| =\|\bv_2 \| = C$
for some positive constant $C$, where
\begin{eqnarray*}
a_n&=&\max_{1\le j\le
p}\{|p_{\lambda_{1j}}'(|\alpha_{0j}|)|,\alpha_{0j}\ne 0\},\\
c_n&=&\max_{1\le k\le q}\{|p_{\lambda_{2k}}'(|\beta_{0k}|)|,\beta_{0k}\ne
0\},
\end{eqnarray*}
and $\alpha_{0j}$ and $\beta_{0k}$ are the $j$th and $k$th elements of
$\balpha_0$ and $\bbeta_0$, respectively, for $j=1,\ldots,p$ and
$k=1,\ldots,q$. Furthermore, define
\begin{eqnarray*}
D_{n,1}&=&\sum_{i=1}^n
\bigl\{Y_i-\wideh\eta\bigl(Z_i^\T(\balpha_0+\gamma_n\bv_1),\balpha_0+\gamma_n\bv_1,
\bbeta_0+\gamma_n\bv_2\bigr)-X_i^\T(\bbeta_0+\gamma_n\bv_2)\bigr\}^2\\
&&{}-\sum_{i=1}^n\{Y_i-\wideh\eta(X_i^\T\balpha_0,\balpha_0,\bbeta_0)-X_i^\T\bbeta_0\}^2
\end{eqnarray*}
and\vspace*{-3pt}
\begin{eqnarray*}
D_{n,2} &=&-n\sum_{j=1}^{p_0}\{p_{\lambda_{1j}}(|\alpha_{0j}+\gamma_n v_{1j}|)-p_{\lambda_{1j}}(|\alpha_{0j}|)\}\\[-2pt]
&&{}-n\sum_{k=1}^{q_0}\{p_{\lambda_{2k}}(|\beta_{0k}+\gamma_n v_{2k}|)-p_{\lambda_{2k}}(|\beta_{0k}|)\}.
\end{eqnarray*}
After algebraic  simplification, we have
%
\begin{eqnarray}\label{A1}
D_{n,1}&=&
\sum_{i=1}^n\bigl\{\wideh\eta\bigl(Z_i^\T(\balpha_0+\gamma_n\bv_1),\balpha_0+\gamma_n\bv_1,\bbeta_0+\gamma_n\bv_2\bigr)\nonumber\\[-2pt]
&&\hphantom{\sum_{i=1}^n}\hspace*{51pt}{}-\wideh\eta(Z_i^\T\balpha_0,\balpha_0,\bbeta_0)+X_i^\T\gamma_n\bv_2\bigr\}\nonumber\\
[-9pt]\\ [-9pt]
&&\hphantom{\sum_{i=1}^n}{}\times\bigl\{\wideh\eta\bigl(Z_i^\T(\balpha_0+\gamma_n\bv_1),\balpha_0+\gamma_n\bv_1,\bbeta_0+\gamma_n\bv_2\bigr)\nonumber\\[-2pt]
&&\hphantom{\sum_{i=1}^n}\hspace*{14pt}{}+\wideh\eta(Z_i^\T\balpha_0,\balpha_0,\bbeta_0)+X_i^\T(\bbeta_0+\gamma_n\bv_2)
+X_i^\T\bbeta_0-2Y_i\bigr\}\nonumber\\[-2pt]
&=&\sum_{i=1}^n\{\wideti Z_i^\T\eta'(\Lambda_i)\bv_1\gamma_n+\wideti X_i^\T\bv_2\gamma_n\}^2\nonumber\\[-2pt]
&&{}-\sum_{i=1}^n\{\wideti
Z_i^\T\eta'(\Lambda_i)\bv_1\gamma_n+\wideti
X_i^\T\bv_2\gamma_n\}\ve_i+o_p(1).\nonumber
\end{eqnarray}
Moreover, applying
the Taylor expansion and the Cauchy--Schwarz inequality, we are able to
show that $n^{-1}D_{n,2}$ is bounded by
\begin{eqnarray*}
&&\sqrt{p_0}\gamma_n a_n
\|\bv_1\|+\gamma_n^2 b_n \|\bv_1\|^2 +\sqrt{q_0} \gamma_n c_n \|\bv_2\|
+\gamma_n^2 d_n\|\bv_2\|^2\\
&&\qquad\le C\gamma_n^2\bigl(\sqrt{p_0}+b_n C+\sqrt{q_0}+d_n C\bigr),
\end{eqnarray*}
where
\[
b_n=\max_{1\le j\le p} \{|p_{\lambda_{1j}}''(|\alpha_{0j}|)|,
\alpha_{0j}\ne 0\}, \qquad d_n=\max_{1\le k\le q}
\{|p_{\lambda_{2k}}''(|\beta_{0k}|)|, \beta_{0k}\ne 0\}.
\]
When $b_n$ and $d_n$ tend to $0$ and $C$ is sufficiently large, the
first term on the right-hand side of (\ref{A1}) dominates the second
term on the right-hand side of~(\ref{A1}) and $D_{n,2}$. As a result,
for any given $\nu > 0$, there exists a large constant $C$ such that
\[
P \Bigl\{\inf_{{\mathcal V}_{12}} {\mathcal L}_P(\balpha_0+\gamma_n \bv_1, \bbeta_0+
\gamma_n \bv_2) > {\mathcal L}_P(\balpha_0,\bbeta_0)\Bigr\} \geq 1 - \nu,
\]
where ${\mathcal V}_{12}=\{(\bv_1,\bv_2)\dvtx\|\bv_1 \| = C, \|\bv_2 \| = C\}$.
We therefore conclude that the rate of convergence of
$(\wideh{\balpha}{}^\T_{\lambda_1}, \wideh{\bbeta}{}^{\T}_{\lambda_2})^{\T}$ is
$O_P(n^{-1/2}+a_n+c_n)$.

\textit{Step II}. Let $\balpha_1$  and  $\bbeta_1$ satisfy
$\|\balpha_1-\balpha_{10} \| = O_P(n^{-1/2})$ and
$\|\bbeta_1-\bbeta_{10} \| = O_P(n^{-1/2})$, respectively. We next show
that
%
\begin{equation} \label{eq:step2}
{\mathcal L}_P\left\{ \pmatrix{\balpha_{1}\cr {\mathbf 0}},
\pmatrix{\bbeta_{1}\cr {\mathbf 0}} \right\} = \min_{\mathcal C}{\mathcal L}_P\left\{\pmatrix{\balpha_{1}\cr \balpha_{2}},
\pmatrix{\bbeta_{1}\cr
\bbeta_{2}}\right\},\vadjust{\goodbreak}
\end{equation}
where ${\mathcal C}=\{\|\balpha_2\|\le C^* n^{-1/2}, \| \bbeta_2 \| \leq
C^* n^{-1/2}\}$ and $C^*$ is a positive constant.\looseness=1

Consider  $\beta_k\in (-C^* n^{-1/2}, C^* n^{-1/2})$ for
$k=q_0+1,\ldots, q$.  When $\beta_k\ne 0$, we have ${\partial {\mathcal
L}_P(\balpha, \bbeta)}/{\partial \beta_k}=\ell_k(\balpha,
\bbeta)+np_{\lambda_{2k}}'(|\beta_k|)\sgn(\beta_k)$, where
\begin{eqnarray*}
\ell_k(\balpha, \bbeta)
&=&\frac1n\sum_{i=1}^n\{Y_i-\wideh\eta(Z_i^\T\balpha,\balpha,\bbeta)+X_i^\T\bbeta\}
    \biggl\{X_{ik}+\frac{\partial\wideh\eta(Z_i^\T\balpha,\balpha,\bbeta)}{\partial\beta_k}\biggr\}\\
&=&\frac1n\sum_{i=1}^n\{\wideh\eta(Z_i^\T\balpha,\balpha,\bbeta)-\eta(Z_i^\T\balpha_0)+
X_i^\T(\bbeta-\bbeta_0)-\ve_i\}\\
&&\phantom{\frac1n\sum_{i=1}^n}{}\times\biggl\{X_{ik}+\frac{\partial\wideh\eta(Z_i^\T\balpha,\balpha,\bbeta)}{\partial\beta_k}\biggr\}.
\end{eqnarray*}
Applying  arguments similar to those used in the proof of Theorem 5.2
of Ichimura~(\citeyear{ich1993}), together with algebraic simplifications, the above
term can be expressed as
\[
\frac1n\sum_{i=1}^n(\balpha-\balpha_0)^\T\wideti Z_i\eta'(\Lambda_i)\wideti
X_{ik} +\frac1n\sum_{i=1}^n(\bbeta-\bbeta_0)^\T\wideti X_i\wideti
X_{ik}+o_p(n^{-1/2}).
\]
Using the assumptions that
$\|\balpha-\balpha_0\|=O_P(n^{-1/2})$ and
$\|\bbeta-\bbeta_0\|=O_P(n^{-1/2})$, we have that
$n^{-1}\ell_k(\balpha, \bbeta)$ is of the order $O_P(n^{-1/2})$.
Therefore,
\[
\frac{\partial \mathcal{L}_P(\balpha,\bbeta)}{\partial \beta_k}= -n
\lambda_{2k}
\{\lambda_{2k}^{-1}p_{\lambda_{2k}}'(|\beta_k|)\sgn(\beta_k) +
O_P(n^{-1/2}\lambda_{2k}^{-1})\}.
\]
Because $\liminf_{n\to\infty}\liminf_{\beta_j\to
0^{+}}\lambda_{2k}^{-1}p_{\lambda_{2k}}'(|\beta_k|)>0$ and
$n^{-1/2}/{\lambda_{2k}}\to 0$,\break ${\partial {\mathcal L}_P(\balpha,
\bbeta)}/{\partial \beta_k}$ and $\beta_k$ have different signs for
$\beta_k\in (-C^* n^{-1/2}, C^* n^{-1/2})$. Analogously, we can show
that ${\partial {\mathcal L}_P(\balpha, \bbeta)}/{\partial \alpha_j}$ and
$\alpha_j$ have different signs when $\alpha_j\in (-C^* n^{-1/2}, C^*
n^{-1/2})$ for $j=p_0+1,\ldots,p$. Consequently, the minimum is
attained  at $\balpha_2=0$ and $\bbeta_2 = 0$. This completes the proof
of (\ref{eq:step2}).

\textit{Step III.} Finally, we demonstrate the asymptotic
normality of $\wideh\balpha_{{\lambda_1}1}$ and
$\wideh\bbeta_{{\lambda_2}1}$. For the sake of simplicity, we define
\begin{eqnarray*}
{\bR_{\lambda_1}}&=&\{p_{\lambda_{11}}'(|\alpha_{01}|)\sgn(\alpha_{01}),\ldots,
p_{\lambda_{1{p_0}}}'(|\alpha_{0{p_0}}|)\sgn(\alpha_{0{p_0}})\}^\T, \\
\Sigma_{\lambda_1}&=&\diag\{p_{\lambda_{11}}''(|\alpha_{01}|),\ldots,p_{\lambda_{1{p_0}}}''(|\alpha_{0{p_0}}|)\},\\
{\bR_{\lambda_2}}&=&\{p_{\lambda_{21}}'(|\beta_{01}|)\sgn(\beta_{01}),\ldots,
p_{\lambda_{2{q_0}}}'(|\beta_{0{q_0}}|)\sgn(\beta_{0{q_0}})\}^\T,\\
\Sigma_{\lambda_2}&=&\diag\{p_{\lambda_{21}}''(|\beta_{01}|),\ldots,p_{\lambda_{2{q_0}}}''(|\beta_{0{q_0}}|)\},
\end{eqnarray*}
where $\alpha_{0j}$ and $\beta_{0k}$ are the $j$th and $k$th
elements of $\balpha_{10}$ and $\bbeta_{10}$, respectively, for
$j=1,\ldots,p_0$ and $k=1,\ldots,q_0$.
It follows from (\ref{eq:2.5}) that $\wideh\balpha_{{\lambda_1}1}$ and
$\wideh\bbeta_{{\lambda_2}1}$ satisfy
%
\begin{eqnarray}\label{eq:A.4}
{\mathbf 0}&=&\left\{
\matrix{\dfrac{\partial {\mathcal L}_P(\wideh\balpha_{{\lambda_1}1},
\wideh\bbeta_{{\lambda_2}1})}{\partial\balpha_1}\vspace*{2pt}\cr
\dfrac{\partial {\mathcal L}_P(\wideh\balpha_{{\lambda_1}1},
\wideh\bbeta_{{\lambda_2}1})}{\partial\bbeta_1}}\right\}\nonumber\\
[-8pt]\\ [-8pt]
&=&l(\wideh\balpha_{{\lambda_1}1},
\wideh\bbeta_{{\lambda_2}1})+\left\{\matrix{
\mathbf{R}_{\lambda_1}-\Sigma_{\lambda_1}(\wideh\balpha_{{\lambda_1}1}-\balpha_{10})\vspace*{2pt}\cr
\mathbf{R}_{\lambda_2}-\Sigma_{\lambda_2}(\wideh\bbeta_{{\lambda_2}1}-\bbeta_{10})}\right\},\nonumber
\end{eqnarray}
where
\[
l(\wideh\balpha_{{\lambda_1}1},
\wideh\bbeta_{{\lambda_2}1})= \frac1n\sum_{i=1}^n
\{Y_i-\wideh\eta(\Lambda_{i},\wideh\bzeta_{\lambda
1})-X_{i,1}^\T\wideh\bbeta_{{\lambda_2}1}\}\left\{\matrix{
\dfrac{\partial \wideh\eta(\wideh\Lambda_{i},\wideh\bzeta_{\lambda 1})}{\partial\balpha}\vspace*{2pt}\cr
\dfrac{\partial \wideh\eta(\wideh\Lambda_{i},\wideh\bzeta_{\lambda 1})}{\partial\bbeta}+X_{i,1}}\right\},
\]
$\wideh\Lambda_{i}$ and $X_{i,1}$ here being the $i$th rows of   $\wideh\Lambda$ and $X_{1}$, respectively,
and $\wideh\bzeta_{\lambda
1}=(\wideh\balpha{}^{\T}_{{\lambda_1}1},\wideh{\bbeta}{}^{\T}_{{\lambda_2}1})^{\T}$
being  the penalized least-squares estimator of
$\bzeta_1=(\balpha_1^{\T},\bbeta_1^{\T})^{\T}$.

Applying the Taylor expansion, we obtain
\begin{eqnarray*}
l(\wideh\balpha_{{\lambda_1}1}, \wideh\bbeta_{{\lambda_2}1})&=&\frac1n\sum_{i=1}^n\{Y_i-\wideh\eta(\Lambda_{i},\bzeta_1)-X_{i,1}^\T\bbeta_{1}\}
\left\{\matrix{
\dfrac{\partial \wideh\eta(\Lambda_{i},\bzeta_1)}{\partial\balpha_1}\vspace*{2pt}\cr
\dfrac{\partial \wideh\eta(\Lambda_{i},\bzeta_1)}{\partial\bbeta_1}+X_{i,1}}\right\}\\
&&{}-\frac1n\sum_{i=1}^n \left\{\matrix{
\dfrac{\partial
\wideh\eta(\bar\Lambda_{i},\bar\bzeta_1)}{\partial\balpha_1}\vspace*{2pt}\cr
\dfrac{\partial \wideh\eta(\bar\Lambda_{i},\bar\bzeta_1)}{\partial\bbeta_1}+X_{i,1}}\right\}^{\otimes2}(\wideh\bzeta_{\lambda 1}-\bzeta_1)\\
&&{}-\frac1n\sum_{i=1}^n \{Y_i-\wideh\eta(\bar\Lambda_{i},
\bar\bzeta_1)-X_{i,1}^\T\bar\bbeta_1\} \frac{\partial^2
\wideh\eta(\bar\Lambda_{i},\bar\bzeta_1)}{\partial\bzeta_1\,\partial\bzeta_1^\T}
(\wideh\bzeta_{\lambda 1}-\bzeta_1),
\end{eqnarray*}
where $\bar\bbeta_1$,
$\bar\Lambda_{i}$  and $\bar\bzeta_1$ are the interior points between
$\bbeta_1$ and $\hat\bbeta_{{\lambda_2}1}$, $\Lambda_{i}$ and
$\hat\Lambda_{i}$, and $\bzeta_1$ and $\hat\bzeta_{\lambda 1}$,
respectively. Furthermore, using   arguments  similar to those used in
the proof of Theorem~\ref{th:th1}, we have that
\begin{eqnarray*}
l(\wideh\balpha_{{\lambda_1}1}, \wideh\bbeta_{{\lambda_2}1}) &=&\frac1{\sqrt
n}\sum_{i=1}^n \left\{\matrix{
\wideti Z_{i,1}\eta'(\Lambda_{i})\vspace*{2pt}\cr
\wideti X_{i,1}}\right\}\ve_i\\
&&{}-\frac1n \sum_{i=1}^n \left\{\matrix{
\wideti Z_{i,1}\eta'(\Lambda_{i})\vspace*{2pt}\cr
\wideti X_{i,1}}\right\}^{\otimes2}\sqrt n\pmatrix{
\wideh\balpha_{{\lambda_1}1}-\balpha_{10}\vspace*{2pt}\cr
\wideh\bbeta_{{\lambda_2}1}-\bbeta_{10}}+o_p(1).
\end{eqnarray*}
Moreover, the summand of the matrix over $n$ in the second term of
the above equation converges to $\left({{\GZZ \enskip \GXZ}\atop {\GXZ^{\T}\enskip
\GXX}}\right)$.
These results, together with  (\ref{eq:A.4}), lead to
\begin{eqnarray*}
&&n^{1/2} \pmatrix{
\GZZ+\Sigma_{\lambda_1} & \GXZ\vspace*{2pt}\cr
\GXZ^{\T}  & \GXX+\Sigma_{\lambda_2}}
\pmatrix{
\wideh\balpha_{{\lambda_1}1}-\balpha_{10}\vspace*{2pt}\cr
\wideh\bbeta_{{\lambda_2}1}-\bbeta_{10}}
-n^{1/2} \pmatrix{
{\bR_{\lambda_1}}\vspace*{2pt}\cr
{\bR_{\lambda_2}}}\\
&&\qquad=\frac1{\sqrt n}\sum_{i=1}^n \left\{\matrix{
\wideti Z_{i,1}\eta'(\Lambda_{i})\ve_i\vspace*{2pt}\cr
\wideti X_{i,1}\ve_i}\right\}+o_p(1).
\end{eqnarray*}
It follows that
\begin{eqnarray*}
&&n^{1/2}(\GZZ+\Sigma_{\lambda_1})(\wideh\balpha_{{\lambda_1}1}-\balpha_{10})
+n^{1/2}\GXZ(\wideh\bbeta_{{\lambda_2}1}-\bbeta_{10})-n^{1/2}{\bR_{\lambda_1}}\\
&&\qquad=\frac1{\sqrt n}\sum_{i=1}^n\wideti Z_{i,1}\eta'(\Lambda_{i})\ve_i+o_p(1)
\end{eqnarray*}
and
\begin{eqnarray*}
&&n^{1/2}\GXZ^{\T}(\wideh\balpha_{{\lambda_1}1}-\balpha_{10})
+n^{1/2}(\GXX+\Sigma_{\lambda_2})(\wideh\bbeta_{{\lambda_2}1}-\bbeta_{10})-n^{1/2}{\bR_{\lambda_2}}\\
&&\qquad= \frac1{\sqrt n}\sum_{i=1}^n\wideti
X_{i,1}\ve_i+o_p(1).
\end{eqnarray*}
After simplification, we have
%
\begin{eqnarray}\label{eq:A.6}\quad
&&\sqrt n\{(\GZZ+\Sigma_{\lambda_1})-\GXZ(\GXX+\Sigma_{\lambda_2})^{-1}\GXZ^\T\}
(\wideh\balpha_{{\lambda_1}1}-\balpha_{10}) \nonumber\\
&&\quad{}+n^{1/2}\{\GXZ(\GXX+\Sigma_{\lambda_1})^{-1}{\bR_{\lambda_1}}-{\bR_{\lambda_2}}\}\\
&&\qquad=\frac1{\sqrt n}\sum_{i=1}^n \{\wideti
Z_{i,1}\eta'(\Lambda_{i})-\GXZ(\GXX+\Sigma_{\lambda_2})^{-1}\wideti
X_{i,1}\}\ve_i+o_p(1)\nonumber
\end{eqnarray}
and
%
\begin{eqnarray}\label{eq:A.7}\quad
&&\sqrt n
\{(\GXX+\Sigma_{\lambda_2})-\GXZ^\T(\GZZ+\Sigma_{\lambda_1})^{-1}\GXZ\}
(\wideh\bbeta_{{\lambda_2}1}-\bbeta_{10})\nonumber\\
&&\quad{}+n^{1/2}\{\GXZ^\T(\GZZ+\Sigma_{\lambda_1})^{-1}{\bR_{\lambda_2}}-{\bR_{\lambda_1}}\}\\
&&\qquad=\frac1{\sqrt n}\sum_{i=1}^n \{\wideti
X_{i,1}-\GXZ^{\T}(\GZZ+\Sigma_{\lambda_1})^{-1}\wideti
Z_{i,1}\eta'(\Lambda_{i})\} \ve_i+o_p(1).\nonumber
\end{eqnarray}
Equations
(\ref{eq:A.6}) and (\ref{eq:A.7}), together with the central limit
theorem, yield that
\begin{eqnarray*}
&&\sqrt n\{\GZZ+\Sigma_{\lambda_1}-\GXZ(\GXX+\Sigma_{\lambda_2})^{-1}\GXZ^\T\}
(\wideh\balpha_{{\lambda_1}1}-\balpha_{10})\\
&&\qquad{}+n^{1/2}\{\GXZ(\GXX+\Sigma_{\lambda_2})^{-1}\bR_{\lambda_1}-\bR_{\lambda_2}\}\to
N(0, \Sba)
\end{eqnarray*}
and
\begin{eqnarray*}
&&\sqrt n\{\GXX+\Sigma_{\lambda_2}-\GXZ^\T(\GZZ+\Sigma_{\lambda_1})^{-1}\GXZ\}
(\wideh\bbeta_{{\lambda_2}1}-\bbeta_{10})\\[-2pt]
&&\qquad{}+n^{1/2}\{\GXZ^\T(\GZZ+\Sigma_{\lambda_1})^{-1}\bR_{\lambda_2}-\bR_{\lambda_1}\}\to
N(0,\Sbb),
\end{eqnarray*}
where
\begin{eqnarray*}
\Sba&=&\{\GZZ-2\GXZ(\GXX+\Sigma_{\lambda_2})^{-1}\GXZ^\T\\[-2pt]
&&{}+\GXZ(\GXX+\Sigma_{\lambda_2})^{-1}\GXX(\GXX+\Sigma_{\lambda_2})^{-1}\GXZ^\T\}\sigma^2
\end{eqnarray*}
and
\begin{eqnarray*}
\Sbb&=&\{\GXX-2\GXZ^\T(\GZZ+\Sigma_{\lambda_1})^{-1}\GXZ\\[-2pt]
&&{}+\GXZ^\T(\GZZ+\Sigma_{\lambda_1})^{-1}\GZZ(\GZZ+\Sigma_{\lambda_1})^{-1}\GXZ\}\sigma^2.
\end{eqnarray*}
Because each element of $n^{1/2}\Sigma_{\lambda_1}$,
$n^{1/2}\Sigma_{\lambda_2}$,  $n^{1/2}\bR_{\lambda_1}$ and
$n^{1/2}\bR_{\lambda_2}$ tends to zero, we complete the proof.

\subsection{\texorpdfstring{Proof of Theorem \protect\ref{th:BIC}}{Proof of Theorem 3}}

Let $\tau_n=\log(n)$, $\lambda_{n1j}=\tau_n\SE(\hat{\alpha}_j^u)$,
 $\lambda_{n2k}=\tau_n\SE(\hat{\beta}_k^u)$ and
\[
\BIC(S_T) = \log\{\sigma_n^2(S_T)\} + \{\log(n)/n\}\df(S_T),
\]
where $\df(S_T)$ stands for the degrees of freedom of the true model
$S_T$. $\SE(\hat{\alpha}_k^u)=O(1/\sqrt{n})$ and
$\SE(\hat{\beta}_j^u)=O(1/\sqrt{n})$. Thus,
$\lambda_{n1j}=O_P(\log(n)/\sqrt{n})$ and
$\lambda_{n2k}=O_P(\log(n)/\sqrt{n})$. Then, employing  techniques
similar to those used  in Wang, Li and Tsai (\citeyear{wanlitsai2007}), we obtain that
%
\begin{equation}\label{eq:bic1}
P\{\BIC(\tau_n)=\BIC(S_T)\}=1.
\end{equation}
Therefore, to prove the theorem, it suffices to show that
%
\begin{equation}\label{eq:bic2}
P\Bigl\{\inf_{\lambda\in \Omega_-\cup \Omega_+}\BIC(\lambda)>\BIC(\tau_n)\Bigr\}
\to 1,
\end{equation}
where
\[
\Omega_-=\{\lambda\dvtx S_\lambda \not\supset S_T\} \quad \mbox{and} \quad
\Omega_+=\{\lambda\dvtx S_\lambda \supset S_T\}
\]
represent the underfitted and overfitted models, respectively.

To demonstrate (\ref{eq:bic2}), we consider two separate cases given
below.

\textit{Case 1: Underfitted model} (i.e., the model misses at least one
covariate from the true model). For any $\lambda\in\Omega_-$,
(\ref{eq:bic1}). together with assumptions~(A) and (B), implies that,
with probability tending to one,
\begin{eqnarray*}
\BIC(\lambda) -\BIC(\tau_n)&=& \log\{\mse(\lambda)\} + \{\log(n)/n\}\df_{\lambda} -\BIC(S_T)\\
&\ge &\log\{\mse(\lambda)\}-\BIC(S_T)\\[-2pt]
&\ge & \log\{\sigma_n^2(S_\lambda)\} -\BIC(S_T)\\
&\ge & \inf_{\lambda\in \Omega_-} \log\{\sigma_n^2(S_\lambda)\}
-\BIC(S_T)\\
&\ge & \min_{S\not\supset S_T} \log\{\sigma_n^2(S)\} -\log\{\sigma_n^2(S_T)\} - \{\log(n)/n\}\df(S_T)\\
&\to & \min_{S\not\supset S_T} \log\{\sigma^2(S)/\sigma^2(S_T)\}\ge 0.
\end{eqnarray*}
\textit{Case 2: Overfitted model} (i.e., the model contains all of the
covariates in the true model and includes at least one covariate that
does not belong to the true model). For any $\lambda\in \Omega_+$, it
follows by (\ref{eq:bic1}) that, with probability tending to one,
\begin{eqnarray*}
n\{\BIC(\lambda)-\BIC(\tau_n)\}&=&n\{\BIC(\lambda)-\BIC(S_T)\}\\
&=& n \log\{\mse(\lambda)/\sigma_n^2(S_T)\} +
(\df_{\lambda}-\df_{S_T})\log(n)\\
&\ge& n \log\{\sigma^2_n(S_\lambda)/\sigma_n^2(S_T)\} +
(\df_{\lambda}-\df_{S_T})\log(n)\\
&=&\frac{n\{\sigma^2_n(S_\lambda)-\sigma^2_n(S_T)\}}{\sigma^2_n(S_T)}\{1+
o_P(1)\} \\
&&{}+ (\df_{\lambda}-\df_{S_T})\log(n).
\end{eqnarray*}
Applying the result of Theorem~\ref{th:test1}, we know that
${n\{\sigma^2_n(S_\lambda)-\sigma^2_n(S_T)\}}/{\sigma^2_n(S_T)}$ is an
asymptotically chi-squared distribution with $\df_{\lambda}-\df_{S_T}$
degrees of freedom. Accordingly, we obtain that
${n\{\sigma^2_n(S_\lambda)-\sigma^2_n(S_T)\}}/{\sigma^2_n(S_T)} =
O_P(1)$. Moreover, for any $\lambda\in \Omega_+$,
$\df_{\lambda}-\df_{S_T}\ge 1$, and hence
$(\df_{\lambda}-\df_{S_T})\log(n)$ diverges to $+\infty$ as
$n\to\infty$. Consequently,
\[
P\Bigl\{\inf_{\lambda\in \Omega_+}
n\{\BIC(\lambda)- \BIC(\tau_n)\} >0\Bigr\} =P\Bigl\{\inf_{\lambda\in \Omega_+}
{\BIC(\lambda)- \BIC(\tau_n)} >0\Bigr\}\to 1.
\]
The results of Cases 1
and 2 complete the proof.

\subsection{\texorpdfstring{Proof of Theorem \protect\ref{th:test1}}{Proof of Theorem 4}}

We apply techniques similar to those used in the proofs of Theorems 3.1
and 3.2 in Fan and Huang (\citeyear{fanhua2005}) to show this theorem. Accordingly, we
only provide a sketch of a proof here;  detailed derivations can be
obtained from the authors upon request.

Let
\[
B_n=\left[\sum_{i=1}^n \left\{\matrix{
\eta'(\Lambda_i)\wideh Z_i\vspace*{2pt}\cr
\wideh X_i}\right\}^{\otimes2}\right]^{-1}
\bA^\T\left(\bA\left[\sum_{i=1}^n \left\{\matrix{
\eta'(\Lambda_i)\wideh Z_i\vspace*{2pt}\cr
\wideh X_i}\right\}^{\otimes2}\right]^{-1}
\bA^\T\right)^{-1}\bA.
\]
The difference $Q(H_0)-Q(H_1)$ can be expressed as
\begin{eqnarray} \label{eq:t1}
&&\sum_{i=1}^n\{Y_i-X_i^\T\wideh\bbeta_0-\wideh\eta(Z_i^\T\wideh\balpha_0,\wideh\bzeta_0)\}^2
 -\sum_{i=1}^n\{Y_i-X_i^\T\wideh\bbeta_1-\wideh\eta(Z_i^\T\wideh\balpha_1,\wideh\bzeta_{1})\}^2 \nonumber \\
&&\qquad=\sum_{i=1}^n\{\wideh\eta(Z_i^\T\wideh\balpha_1,\wideh\bzeta_{1})-\wideh\eta(Z_i^\T\wideh\balpha_0,
\wideh\bzeta_0)+X_i^\T(\wideh\bbeta_1-\wideh\bbeta_0)\}^2\nonumber\\[-2pt]
&&\qquad\quad{}+2\sum_{i=1}^n\{Y_i-X_i^\T\wideh\bbeta_1-\wideh\eta(Z_i^\T\wideh\balpha_1,\wideh\bzeta_{1})\}\\[-2pt]
&&\qquad\quad\phantom{{}+2\sum_{i=1}^n}{}\times\{\wideh\eta_1(Z_i^\T\wideh\balpha_1,\wideh\bzeta_{1})-\wideh\eta(Z_i^\T\wideh\balpha_0, \wideh\bzeta_0)
+X_i^\T(\wideh\bbeta_1-\wideh\bbeta_0)\}\nonumber\\[-2pt]
&&\qquad\stackrel{\operatorname{def}}{=} Q_1+Q_2.\nonumber
\end{eqnarray}
It can be shown that  $Q_2$ is asymptotically
negligible in probability. Furthermore, $Q_1$ can be simplified as
\[
Q_1 =(\wideh\bzeta_{1}-\wideh\bzeta_0)^\T \sum_{i=1}^n \left\{
\matrix{
\eta'(Z_i^\T\wideh\balpha_0,\wideh\bzeta_0)\wideh Z_i\vspace*{2pt}\cr
\wideh X_i}\right\}^{\otimes2}
(\wideh\bzeta_{1}-\wideh\bzeta_0)+o_p(1).
\]
A direct calculation yields
that $\wideh\bzeta_0=\wideh\bzeta_{1}-B_n\wideh\bzeta_{1}$. This, together with
the asymptotic normality and consistency of $\hat{\bzeta}_1$ obtained
from Theorem~\ref{th:th1}, implies that $\sigma^{-2}Q_1\to \chi_m^2$ in
distribution under $H_0$. Moreover, under $H_1$, $\sigma^{-2}Q_1$
asymptotically follows a noncentral chi-squared distribution with $m$
degrees of freedom and noncentrality parameter $\phi$. This completes
the proof.

\subsection{\texorpdfstring{Proof of Theorem \protect\ref{th:glrt}}{Proof of Theorem 5}}

It is noteworthy that
\begin{eqnarray*}
\rss(H_0)&=& \sum_{i=1}^n
    \{Y_i-X_i^\T\bbeta-\wideti\eta(Z_i^\T\wideh\balpha)\}^2\\[-2pt]
&&{}+\sum_{i=1}^n
    [\{Y_i-X_i^\T\wideh\bbeta- \wideh\eta(Z_i^\T\wideh\balpha)\}^2
        -\{Y_i-X_i^\T\bbeta-\wideti\eta(Z_i^\T\wideh\balpha)\}^2]\\[-2pt]
&\stackrel{\operatorname{def}}{=}& \rss^*(H_0)+I_{n0}
\end{eqnarray*}
and
\begin{eqnarray*}
\rss(H_1)&=& \sum_{i=1}^n
    \{Y_i-X_i^\T\bbeta-\wideh\eta(Z_i^\T\wideh\balpha)\}^2\\[-2pt]
&&{}+\sum_{i=1}^n
    [\{Y_i-X_i^\T\wideh\bbeta- \wideh\eta(Z_i^\T\wideh\balpha)\}^2
    -\{Y_i-X_i^\T\bbeta-\wideh\eta(Z_i^\T\wideh\balpha)\}^2]\\[-2pt]
&\stackrel{\operatorname{def}}{=} &\rss^*(H_1)+I_{n1}.
\end{eqnarray*}
Applying similar arguments to those
in the proof of Theorem 5 in Fan, Zhang and Zhang (\citeyear{fanzhazha2001}), under $H_0$,
we have
\[
\frac{n r_K}{2} \frac{\rss^*(H_0)-\rss^*(H_1) } {\rss^*(H_1)}
\sim \chi_{df_n}^2,\vadjust{\goodbreak}
\]
where $df_n$ is defined as in
Theorem~\ref{th:glrt} and it approaches infinity as $n\to \infty$.
Furthermore, it can be  straightforwardly shown that
$n^{-1}I_{n0}=\sigma^2\{1+o_P(1)\}$ and
$n^{-1}I_{n1}=\sigma^2\{1+o_P(1)\}$. These results complete the proof.
\end{appendix}

\section*{Acknowledgments}
The authors wish to thank the former Editor, Professor Susan Murphy,
an Associate Editor and three referees for their constructive comments
that substantially improved an earlier version of this paper.

\printaddresses

\end{document}